\ifx\LocalElevenPtArticleFormatLoaded\undefined
\documentclass[a4paper,11pt]{article}     
\usepackage{amsmath,amscd,amssymb}        
\usepackage{latexsym}
\usepackage[german,english]{babel}
\fi
\usepackage[applemac]{inputenc}
\usepackage{epsfig}
\usepackage{xypic}

\usepackage{theorem}

\newtheorem{lemma}{Lemma}[section]

\newtheorem{proposition}[lemma]{Proposition}
\newtheorem{theorem}[lemma]{Theorem}
\newtheorem{cor}[lemma]{Corollary}

\theorembodyfont{\rm}
\newtheorem{remark}[lemma]{Remark}

\newcommand{\op}{\operatorname}
\newenvironment{Proof}{\begin{ProofwCaption}{Proof}}{\end{ProofwCaption}}
\newenvironment{ProofwCaption}[1]%
  {\addvspace\theorempreskipamount \noindent{\it #1.}\rm}
  {\qed \par \addvspace\theorempostskipamount}

\newcommand{\CC}{\mathbb C}

\newcommand{\HH}{\mathbb H}

\newcommand{\LL}{\mathbb L}
\newcommand{\MM}{\mathbb M}
\newcommand{\NN}{\mathbb N}

\newcommand{\PP}{\mathbb P}
\newcommand{\QQ}{\mathbb Q}
\newcommand{\RR}{\mathbb R}

\newcommand{\ZZ}{\mathbb Z}
\newcommand{\cA}{\mathcal A}

\renewcommand{\cD}{\mathcal D}

\newcommand{\cF}{\mathcal F}

\newcommand{\cM}{\mathcal M}

\newcommand{\cO}{\mathcal O}
\newcommand{\cP}{\mathcal P}

\newcommand{\cS}{\mathcal S}

\newcommand{\To}{\longrightarrow}

\newcommand{\Sum}{\sum\limits}

\newcommand{\tens}{\otimes}

\renewcommand{\Bar}{\overline}
\newcommand{\hf}{{\sfrac{1}{2}}}
\newcommand{\half}{{\frac{1}{2}}}

\newcommand{\imic}{\cong}
\newcommand{\GL}{\mathop{\mathrm {GL}}\nolimits}

\renewcommand{\Sp}{\mathop{\mathrm {Sp}}\nolimits}

\newcommand{\Nef}{\mathop{\mathrm {Nef}}\nolimits}
\newcommand{\Pic}{\mathop{\mathrm {Pic}}\nolimits}
\newcommand{\Proj}{\mathop{\mathrm {Proj}}\nolimits}

\newcommand{\Supp}{\mathop{\mathrm {Supp}}\nolimits}
\newcommand{\Sym}{\mathop{\mathrm {Sym}}\nolimits}

\newcommand{\Hom}{\mathop{\mathrm {Hom}}\nolimits}
\newcommand{\Ker}{\mathop{\mathrm {Ker}}\nolimits}

\newcommand{\depth}{\mathop{\mathrm {depth}}\nolimits}

\newcommand{\id}{\mathop{\mathrm {id}}\nolimits}

\newcommand{\pr}{\mathop{\mathrm {pr}}\nolimits}

\newcommand{\Tr}{\mathop{\mathrm {Tr}}\nolimits}

\newcommand{\bde}{\mathbf e}

\newcommand{\bds}{\mathbf s}

\newcommand{\bdx}{\mathbf x}
\newcommand{\bdy}{\mathbf y}

\newcommand{\gothV}{\mathfrak V}

\newcommand{\gothS}{\mathfrak S}

\newcommand{\ratmap}{\dasharrow}
\newcommand{\Eins}{{\mathbf 1}}

\newcommand{\Vor}{{\op{Vor}}}
\newcommand{\Igu}{{\op{Igu}}}
\newcommand{\Sat}{{\op{Sat}}}
\newcommand{\Star}{{\op{Star}}}
\newcommand{\Span}[1]{{\left\langle{#1}\right\rangle}}
\newcommand{\AVOR}[1]{{\cA_{#1}^{\Vor}}}
\newcommand{\AIGU}[1]{{\cA_{#1}^{\Igu}}}
\newcommand{\ASAT}[1]{{\cA_{#1}^{\Sat}}}
\newcommand{\DVOR}[1]{{D_{#1}^{\Vor}}}
\newcommand{\DIGU}[1]{{D_{#1}^{\Igu}}}
\newcommand{\Psing}{{P_{\op{sing}}}}
\newcommand{\RRp}{{\RR_{\ge 0}}}
\newcommand{\emin}{{e_{\min}}}
\newcommand{\Temb}{{\mathrm {emb}\,}}
\newcommand{\sfrac}[2]{{\textstyle{\frac{#1}{#2}}}}
\newcommand{\pxi}{{p_{\xi}^*}}

\newcommand{\qedsymbol}{\mbox{$\Box$}}
\newcommand{\qed}{\unskip\nobreak\hfil\penalty50\hskip1em\hbox{}\nobreak
\hfill\qedsymbol\parfillskip=0pt\finalhyphendemerits=0}
\begin{document}

\title{The nef cone of toroidal compactifications of $\cA_4$}
\author{K.~Hulek and G.K.~Sankaran}
\maketitle

\section*{0 Introduction}
The moduli space $\cA_g$ of principally polarised abelian $g$-folds is
a quasi-projective variety. It has a natural projective compactification, the
Satake compactification, which has bad singularities at infinity. By the
method of toroidal compactification we can construct other
compactifications with milder singularities, at the cost of some loss
of uniqueness. Two popular choices of toroidal compactification are
the Igusa and the Voronoi compactifications: these agree for $g\le 3$
but for $g=4$ they are different.

In this paper, we shall be mainly interested in the Voronoi
compactification $\AVOR4$ of $\cA_4$. This is a natural choice from
the point of view of moduli in view of the results of Alexeev and
Nakamura~(\cite{Al}, \cite{AN}), who show that $\AVOR{g}$ represents a
functor of geometric interest. The case $g=4$ is also of particular
interest as it is the first case where the Torelli map is not dominant
and where we therefore cannot use results from the moduli space of
curves.

In our main result, Theorem~\ref{mainthm}, we
describe the cones of nef divisors on $\AIGU4$ and $\AVOR4$. The proofs are
inductive in the sense that they involve a reduction to the cases
$g=3$ and $g=2$, where comparable results already exist; but some new
techniques are also necessary for the proof. 

However, the
Voronoi compactification for $g>4$ is rather complicated and for this
reason we are not at present able to extend our results even to
$g=5$. We also show (Theorem~\ref{Kample}) that the canonical bundle
on $\AIGU4(n)$ is ample for~$n\ge 3$.

The paper is structured as follows. Section~\ref{toroidal} covers the
facts we need about the different toroidal compactifications that are
available. We describe the Voronoi compactification, in particular, in
some detail, and state the main results. In
Section~\ref{nefconeofpartial} we explain what is known about the
partial compactification of Mumford, which we shall need later. In
Section~\ref{Vorboundary} we describe the fine structure of the
Voronoi boundary in the case $g=4$, which is largely a matter of
understanding the behaviour over the lowest stratum of the Satake
compactification $\ASAT4$. The methods here are toric and much is
deduced from the combinatorics of a single cone in a certain
$10$-dimensional real vector space. The main technical result is that
each non-exceptional boundary divisor of $\AVOR4(n)$, where $n\ge 3$
is a level structure, has a fibration over $\AVOR3(n)$. This is the
inductive step that allows us to deduce facts about $\AVOR4$ from the
cases where~$g<4$. Finally, in Section~\ref{mainsection}, we assign to
a curve in $\AVOR4$ an invariant called the depth, which is the
stratum of $\ASAT4$ that it comes from, and work through the five
cases $0\le \depth(C)\le 4$ that arise. No two of the cases turn out
to be exactly alike.

\noindent{\it Acknowledgements:\/ }We are grateful
to the DAAD and the British Council for financial assistance under ARC
Project 313-ARC-XIII-99/45; to the research network EAGER, supported
by the programme {\it Improving Human Potential \& the
Socio-economic Knowledge base\/} of the European Commission (Contract
No. HPRN-CT-2000-00099); to the DFG Schwerpunktprogramm ``Globale
Methoden in der komplexe Geometrie'' (grant Hu~337/5-1); and to the
Isaac Newton Institute for Mathematical Sciences in Cambridge. We
benefited from discussions with V.~Alexeev, and E.~Schellhammer and
G.~Starke provided precursors of a computer program we used.

\section{Toroidal compactifications}\label{toroidal}
The moduli space of principally polarised abelian varieties of dimension
$g$ is given (over the complex numbers $\CC$) as the quotient
$$
\cA_g=\Sp (2g, \ZZ)\backslash \HH_g.
$$
We shall also consider full (symplectic) level-$n$ structures. The
corresponding moduli spaces are
$$
\cA_g(n)=\Gamma_g(n)\backslash\HH_g
$$
where $\Gamma_g(n)$ is the principal congruence subgroup of level $n$, i.e.\ 
the set of all matrices $\gamma\in\Sp (2g, \ZZ)$ that are
congruent to the unit matrix $\Eins_{2g}$ mod~$n$.
The varieties $\cA_g(n)$ are quasi-projective, but not projective,
varieties with at most finite quotient singularities. The Satake
compactification $\ASAT{g}$ is the minimal compactification
of $\cA_g$. It is $\Proj$ of the ring of modular forms for $\Sp(2g,\ZZ)$.
Set-theoretically $\ASAT{g}$ is the disjoint union
$$
\ASAT{g}=\cA_g\amalg \cA_{g-1}\amalg\ldots\amalg\cA_0
$$
where $\cA_0$ is a point.

Mumford \cite{Mu} introduced a partial compactification
$$
\cA'_g=\cA_g\amalg D'_g
$$
by adding the rank $1$ degenerations. This is again a quasi-projective,
but not projective, variety. There are several toroidal compactifications
$\cA_g^{\Sigma}$. These depend on the choice of a fan $\Sigma$ in
the cone of positive definite $g\times g$ matrices (see below,
Remark~\ref{fans}, for a more precise explanation).
All of them contain $\cA'_g$. The most important choices are:
\begin{itemize}
\item The perfect cone (or first Voronoi) decomposition: see~\cite{V1}
or, for instance, \cite{Co} for definitions and details.
\item The central cone decomposition~$\Igu(g)$. This leads to the Igusa
compactification $\AIGU{g}$.
\item The second Voronoi decomposition~$\Vor(g)$, defined
in~\cite{V2b}. This leads to the Voronoi compactification $\AVOR{g}$.
\end{itemize}

For $g\le 3$ all these fans coincide. For $g=4$ the perfect cone and
central cone decompositions coincide, but the second Voronoi
decomposition is a refinement of the first two: this means that there
is a birational morphism $\AVOR4\to\AIGU4$. For general $g$ very
little is known explicitly about the decompositions and their relation to
each other. There is always a morphism $\cA_g^\Sigma\to
\ASAT{g}$ for any toroidal compactification, and $\cA'_g$
is the inverse image of
$\cA_g\amalg\cA_{g-1}$ under this morphism.

For $g\le 4$ the above decompositions are explicitly known (see e.g.
\cite{V2b}, \cite{ER1}, \cite{ER2}). Since the fan $\Vor(4)$ is
basic the space $\AVOR4$ has only finite quotient
singularities and $\AVOR4(n)$ is smooth for $n\ge 3$. The
spaces $\AIGU4(n)$ are always singular. We shall denote by
$\DVOR{g}$ and $\DIGU{g}$ the closures of $D'_g$ in
$\AVOR{g}$ and $\AIGU{g}$ respectively. It is well known that
$\AIGU{g}$ is a blow-up of the Satake compactification $\ASAT{g}$ and
Alexeev \cite{Al} has proved the same for $\AVOR{g}$. In particular,
$\DVOR{g}$ and $\DIGU{g}$ are $\QQ\,$-Cartier divisors.
In any case it is clear that $\DVOR4$ is $\QQ\,$-Cartier, since $\AVOR4$
is an orbifold and thus $\QQ\,$-factorial. We can see directly that
$\DIGU4$ is $\QQ\,$-Cartier by exhibiting a suitable support function:
see Remark~\ref{support} below.

We denote by $L$ the $\QQ\,$-line bundle of modular forms of weight~$1$
on $\ASAT{g}$, and also its pullback to $\AVOR{g}$ or to $\AIGU{g}$.

\begin{proposition}\label{PicA'}
$\Pic \cA'_g\otimes \QQ=\QQ D'_g\oplus \QQ L$ for $g\ge 2$.
\end{proposition}

\begin{Proof}
This is proved by Mumford (\cite[p. 355]{Mu}) for $g\ge 4$. It is also
well known for $g=2$ and $g=3$: see for instance~\cite{vdG}.
\end{Proof}

For what follows we shall need explicit descriptions of the perfect cone
(=central cone) and the second Voronoi decompositions in the case~$g=4$.

We fix generators $x_1,\ldots,x_4$ for a free abelian group
$\LL_4\imic\ZZ^4$, and we denote $\Sym_2(\LL_4)$ by $\MM_4$; so
$\MM_4\imic\ZZ^{10}$ is the space of $4\times 4$ integer symmetric matrices
with respect to the basis $x_i$. We shall use the basis of $\MM_4$
given by the matrices $U^\ast_{ij}$, $1\le i\le j\le 4$ given by
$$
(U^\ast_{ij})_{kl}=\delta_{\{i,j\},\{k,l\}}.
$$
Thus $U^\ast_{ii}$ is the diagonal matrix with~$1$ in the $i$th place,
corresponding to the quadratic form~$x_i^2$, and $U^\ast_{ij}$ has~$1$
in the $ij$- and $ji$-places, corresponding to the
quadratic form $2x_ix_j$ for $1\le i<j\le 4$. 

The cone $\Sym_2^+(\LL_4\tens\RR)$ is defined to be the convex hull
(that is, $\RRp$-span) of the positive semidefinite forms in
$\MM_4\tens\QQ$. The perfect cone decomposition and the second Voronoi
decomposition are decompositions of the cone $\Sym_2^+(\LL_4\tens\RR)\subset
\MM_4\tens\RR$ into rational polyhedral cones; that is, polyhedral
cones with generators in~$\MM_4$. These cones form fans $\Igu(4)$
(coming from the perfect cone decomposition) and $\Vor(4)$, which are
invariant under the action of $\GL(\LL_4)\imic\GL(4,\ZZ)$.

\begin{remark}\label{fans}
$\Sym_2^+(\LL_4\tens\RR)$, or more generally
$\Sym_2^+(\LL_g\tens\RR)$, is defined in terms of the
lattice~$\LL_g$, and does not depend just on the vector
space~$\LL_g\tens\RR$. The same is true of the torus embeddings
$T_{\MM_g}\Temb(\Sigma)$ (see~\cite{Oda}) which are defined by fans
$\Sigma$ in
$\Sym_2^+(\LL_g\tens\RR)$ and which are used to construct the
compactifications $\cA_g^\Sigma$. If there is no danger of confusion
about which lattice (and hence which torus) is intended, we sometimes
denote $T_{\MM_g}\Temb(\Sigma)$ by~$X_\Sigma$.
\end{remark}

In any real vector space $V$ (usually $V=\MM_4\tens\RR$ or its dual)
we denote the closed cone
$\RRp q_1+\cdots+\RRp q_k$ generated by
$\{q_1,\ldots,q_k\}\subset V$ by $\Span{q_1,\ldots, q_k}$.
In particular $\Span{\pm q}$ is the line~$\RR\, q$.

The perfect cone decomposition has, up to $\GL(\LL_4)$-equivalence, two
maximal, i.e.\ $10$-dimensional, cones: the principal cone $\Pi_1(4)$ and the
second perfect cone $\Pi_2(4)$. The principal cone is given by
$$
\Pi_1(4)=\Span{x^2_1,\ldots,x^2_4,(x_1-x_2)^2,\ldots,(x_3-x_4)^2}.
$$
This cone is basic. The second perfect cone is given by
\begin{eqnarray*}
\Pi_2(4)&=&
\left\langle
 x_1^2, x_2^2, x_3^2, x_4^2,(x_1-x_3)^2,(x_1-x_4)^2 ,(x_2-x_3)^2,(x_2-x_4)^2,
\right.\\
&&\hskip-20pt\left.(x_3-x_4)^2,(x_1+x_2-x_3)^2,(x_1+x_2-x_4)^2,
(x_1+x_2-x_3-x_4)^2\right\rangle.
\end{eqnarray*}
$\Pi_2(4)$ has $64$ $9$-dimensional faces, which fall into two
$\GL(\LL_4)$-equivalence classes called BF and RT: see
\cite{ER2} and the proof of Proposition~\ref{orbitsofexccones}, below.
Representatives of the orbits are given by setting the coefficients of
$x^2_1$, $x^2_3$ and $x^2_4$, respectively of $(x_2-x_3)^2$,
$(x_2-x_4)^2$ and $(x_1+x_2-x_3-x_4)^2$, equal to~$0$. These cones are
basic. Hence $\AIGU4$ has exactly one singular point, which we
denote~$\Psing$.

In order to describe the second Voronoi decomposition we have to introduce
another ray $\eta$, generated by the sum of the primitive
generators of~$\Pi_2(4)$. The primitive generator $e$ of $\eta$ in~$\MM_4$
is given by
\begin{eqnarray}\label{definee}
e&=&\frac{1}{3}\Big[x^2_1+x^2_2+x^2_3+x^2_4\nonumber\\
&&\mbox{}+(x_1-x_3)^2+(x_1-x_4)^2+(x_2-x_3)^2+(x_2-x_4)^2+(x_3-x_4)^2
\nonumber\\
&&\mbox{}+(x_1+x_2-x_3)^2+(x_1+x_2-x_4)^2+(x_1+x_2-x_3-x_4)^2\Big]
\nonumber\\
&=&2(x^2_1+x^2_2+x^2_3+x^2_4+x_1x_2-x_1x_3-x_1x_4-x_2x_3-x_2x_4).
\end{eqnarray}

The second Voronoi decomposition of $\Sym^+_2(\LL_4\tens\RR)$ is the
refinement of the central cone decomposition given by adding all cones
which arise as the span of the central ray $\eta$ with the
$9$-dimensional faces of $\Pi_2(4)$ and the faces of these cones,
together with their $\GL(\LL_4)$-translates. Up
to $\GL(\LL_4)$ this defines two new $10$-dimensional cones, both of
which are basic. Hence $\AVOR4$ is an orbifold, and there is a map
$\pi\colon\AVOR4\to \AIGU4$ given by blowing up a certain ideal sheaf
$\gothV$ supported at the singular point $\Psing\in\AIGU4$.
Let $E$ be the exceptional divisor of this blow-up, i.e.\ the divisor
corresponding to the ray~$\eta$. Actually $\gothV$ is the maximal ideal
of $\cO_{\AIGU4,\Psing}$ and the singularity at $\Psing$ is the cone
on $E$, but we do not need this fact. It can be deduced, for instance,
from~\cite[Theorem I.10]{TE}.

To simplify some calculations it is also useful to consider the Voronoi
transformation $\Psi\colon\LL_4\to \LL_4$, defined by
\begin{equation}\label{Voronoitransformation}
\Psi\colon (x_1, x_2, x_3, x_4)\longmapsto
(x_1+x_2, x_1-x_2, x_1-x_3, x_1-x_4)
\end{equation}
and the induced embedding
$$
\Psi'=\Sym_2(\Psi)\colon\MM_4=\Sym_2(\LL_4)\To \MM_4.
$$
Note that $\Psi$ and $\Psi'$ are embeddings but not isomorphisms,
since $\det\Psi=2$.
We have
$$
\Psi'\big(\Pi_2(4)\big)=\Span{\left\{(x_i\pm x_j)^2,\ 1\le i<j\le 4\right\}},
$$
so if we put $\bdy=\Psi_\QQ^{-1}(\bdx)$ we may express $\Pi_2(4)$ in the
convenient form
\begin{eqnarray}\label{betarho}
\Pi_2(4)&=&\Span{(y_i\pm y_j)^2, 1\le i< j \le 4}\nonumber\\
&=&\Big\{\Sum_{1\le i <j\le 4} 
\left(\beta_{ij}(y_i+y_j)^2 +\rho_{ij}(y_i-y_j)^2\right) \Big|
\beta_{ij},\rho_{ij}\in\RRp\Big\}.
\end{eqnarray}
The generator $e$ of $\eta$ is mapped to
$$
\Psi'(e)=2(x_1^2+x_2^2+x_3^2+x_4^2).
$$
Now let $\pi\colon\AVOR4(n)\to\AIGU4(n)$ and let
$E(n)$ be the exceptional divisor in $\AVOR4(n)$. We set
$D_4(n)=\pi^* \big(\DIGU4(n)\big)$.
\begin{proposition}\label{D4versusD4Vor}
$D_4(n)=\pi^* \big(\DIGU4(n)\big)=\DVOR4(n)+4E(n)$.
\end{proposition}
\begin{Proof}
The level structure plays no part here so we suppress it, taking $n=1$
without loss of generality and writing $E$ for $E(1)$ and so on.
We shall first consider the toric situation. Let
$$
\Tr\colon\MM_4\imic\Sym_2(\ZZ^4)\to \ZZ
$$
be the linear form given by the trace. Then $\Tr'=\Tr\circ\Psi'$ is an
integral linear
form on $\MM_4$ which is $2$-divisible. The form $\half\Tr'$
assumes the value~$1$ on all basic generators of the $1$-dimensional rays of
$\Pi_2(4)$ and the value~$4$ on the $\MM_4$-primitive generator~$e$
of~$\eta$.

Locally (analytically) near the singular point $P_{\op{sing}}$ on
$\AIGU4$ and near the exceptional locus $E$ in
$\AVOR4$, the moduli spaces $\AIGU4$ and $\AVOR4$ are isomorphic to finite
quotients of the toric varieties $X_{\Igu(4)}$ and $X_{\Vor(4)}$
respectively. The finite group by which we take the quotient is the
stabiliser of $P_{\op{sing}}$, respectively $E$. It is a subgroup of
$\GL(\LL_4)$, which acts on $\Sym^+_2(\LL_4\tens\RR)$
by $M \mapsto {}^tQ^{-1}MQ^{-1}$. It is enough to compute the subgroup which
fixes $E$ pointwise. A straightforward calculation shows that this is
$\pm\Eins_4$, which acts trivially. Together with the above toric
calculation this shows that $\pi^* (\DIGU4)=\DVOR4+4E$.
\end{Proof}
\begin{remark}\label{support}
Considering $\half\Tr'$ as a support function on the fan $\Igu(4)$
shows that the boundary $\DIGU4$ of
$\AIGU4$ is a $\QQ\,$-Cartier divisor and that the boundary of
$\AIGU4(n)$ for $n\ge 3$ is a Cartier divisor.
\end{remark}
\begin{cor}\label{canonicalGorenstein}
Let $n\ge 3$. Then $\AIGU4(n)$ is a Gorenstein variety with
canonical singularities.
\end{cor}
\begin{Proof}
For $n\ge 3$ the group $\Gamma_g(n)$ is neat. Hence we only have to
consider singularities which come from the toric construction. The
varieties $\AIGU4(n)$, $n\ge 3$ are normal varieties with
finitely many singularities. Outside these singularities the canonical
divisor is given by
$$
K=\big(5L-\DIGU4(n)\big)|_{\AIGU{4,\op{smooth}}(n)}
$$
where $\DIGU4(n)$ is the boundary. Both $L$ and $\DIGU4(n)$
are Cartier divisors on $\AIGU4(n)$ and hence
\begin{equation}\label{KAIGU}
K_{\AIGU4(n)}=i_*K=5L-\DIGU4(n)
\end{equation}
where $i$ is the inclusion. This shows that these varieties are
Gorenstein.

The varieties $\AVOR4(n)$, $n\ge 3$ are smooth and the
canonical divisor is
\begin{equation}\label{KAVOR}
K_{\AVOR4(n)}=5L-\DVOR4(n)-\sum_s E_s(n),
\end{equation}
where the $E_s(n)$ are the irreducible exceptional divisors of the blow-up map
$\pi\colon\AVOR4(n)\to\AIGU4(n)$. Since
$$
\pi^*\Big(K_{\AIGU4(n)}\Big)=5L-\DVOR4(n)-\sum_s 4E_s(n)
$$
it follows that $\AIGU4(n)$ has canonical, in fact terminal, singularities.
\end{Proof}

We define the open set $\cA^0_4=\AIGU4\setminus P_{\op{sing}}
=\AVOR4\setminus E$, common to both toroidal compactifications.
\begin{proposition}\label{picard}
The Picard groups satisfy
\begin{eqnarray*}
\Pic\AIGU4\otimes\QQ&\cong&\Pic\cA^0_4\otimes\QQ\ \cong
\ \QQ\,L\oplus\QQ\,\DIGU4,\\
\Pic\AVOR4\otimes\QQ&=&\QQ\,L\oplus\QQ\,\DVOR4\oplus\QQ\,E.
\end{eqnarray*}
\end{proposition}

\begin{Proof}
Restricting line bundles defines maps
$$
\Pic\AIGU4\To\Pic\cA^0_4\To\Pic\cA'_4.
$$
All the varieties involved are normal and since the codimensions of
$\cA^0_4\setminus\cA'_4$ in $\cA^0_4$ and of $\AIGU4\setminus\cA^0_4$
in $\AIGU4$ are at least~$2$, these maps are injective. Since
$\Pic\cA'_4=\QQ L\oplus\QQ D'_4$ and since both $L$ and $D'_4$ extend to
$\QQ\,$-line bundles on $\AIGU4$ these maps are also surjective.

The exceptional locus $E$ is irreducible, being the image of the closure
of a torus orbit. Hence the claim about $\Pic\AVOR4\tens\QQ$ follows from
the exact sequence of Chow groups
$$
A_9(E)\otimes\QQ\to A_9(\AVOR4)\otimes\QQ\to
A_9(\cA^0_4)\otimes\QQ\to 0
$$
(see~\cite[Proposition 1.8]{Ful1}).
\end{Proof}

We are now in a position to state the main results of this paper. The
first result is auxiliary and can be stated for general $g\ge 2$.
Note that although $\cA'_g$ is not a projective variety we can still
speak about nef line bundles. By this we mean line bundles whose
restriction to each complete curve has non-negative degree.
\begin{proposition}\label{nefA'}
The nef cone of $\cA'_g$ for $g\ge 2$ is given by
$$
\Nef(\cA'_g)=\left\{aL-bD'_g\mid\ b\ge 0,\ a\ge 12 b\right\}.
$$
\end{proposition}
For the projective varieties $\AIGU4$ and $\AVOR4$ we obtain much better
results.
\begin{theorem}\label{mainthm}
The nef cone of $\AIGU4$ is given by
$$
\Nef(\AIGU4)=\left\{aL-b\DIGU4\mid b\ge 0,\ a\ge 12b\right\}.
$$
The nef cone of $\AVOR4$ is given by
$$
\Nef(\AVOR4)=\left\{aL-bD_4-cE\mid a\ge 12b,\ b\ge 2c\ge 0\right\}.
$$
\end{theorem}
\begin{remark}\label{nefbyVor}
If we work with $\DVOR4$ rather than $D_4$ then, in view of
Proposition~\ref{D4versusD4Vor}, the nef cone
of the Voronoi compactification has the following description:
$$
\Nef(\AVOR4)=\left\{\alpha L- \beta \DVOR4- \gamma E\mid
\beta \ge 0,\ \alpha \ge 12\beta,\ \gamma \ge 4\beta \ge
{\sfrac{8}{9}} \gamma\right\}.
$$
\end{remark}

\begin{remark}\label{Galois}
We have Galois coverings
$$
\alpha_{n,\Igu}\colon\AIGU4(n)\to\AIGU4,\qquad
\alpha_{n,\Vor}\colon\AVOR4(n)\to\AVOR4
$$
given by an action of $\Sp(8,\ZZ/n)$. These coverings, which extend
the obvious covering $\cA_4(n)\to\cA_4$, exist because the definitions
of perfect cone and Voronoi decomposition (\cite{V1}, \cite{V2b},
\cite{Co}, \cite{ER1}) are purely lattice-theoretic and so the
collections of fans that define the Igusa and Voronoi
compactifications are $\Sp(8,\ZZ)$-invariant. Compare
\cite[Proposition 5.1]{San} for a similar situation in the $g=2$ case.

The inverse images $\DIGU4(n)$ and $\DVOR4(n)$ of $\DIGU4$ and
$\DVOR4$ will have several components, as will the inverse image
$E(n)$ of~$E$. The above Galois covers are ramified of order~$n$ along the
boundary, i.e.\ $\alpha_{n,\Igu}^*\big(\DIGU4\big)=n\DIGU4(n)$
and $\alpha_{n,\Vor}^*\big(\DVOR4\big)=n\DVOR4(n)$; it then follows
from Proposition~\ref{D4versusD4Vor} that $\alpha_{n,\Vor}^*(E)=nE(n)$.
\end{remark}
The Picard groups of $\AIGU4(n)$ and
$\AVOR4(n)$ will be much bigger than those of $\AIGU4$
and $\AVOR4$, but we still obtain a description of part of the nef cone.
\begin{cor}\label{nefatleveln}
A divisor $aL-b\DIGU4(n)$ on $\AIGU4(n)$ is nef if and only
if $b\ge 0$ and $a\ge 12 b/n$.

A divisor $aL-b D_4(n)-cE(n)$ on $\AVOR4(n)$
is nef if and only if $a\ge 12 b/n$ and $b\ge 2c\ge 0$.
\end{cor}
This also allows us to draw a conclusion about the nefness of the canonical
divisor.
\begin{lemma}\label{Katleveln}
For any $n\in\NN$
\begin{eqnarray*}
K_{\AIGU4(n)} &=& 5L-\DIGU4(n),\\
K_{\AVOR4(n)} &=& 5L-\DVOR4(n)-E(n)=5L-D_4(n)+3E(n).
\end{eqnarray*}
\end{lemma}
\begin{Proof}
For $n\ge 3$ this was shown above (equations~\eqref{KAIGU}
and~\eqref{KAVOR}). To show that these equalities also hold for $n=1$
and $n=2$, it is enough to check that there are no elements in
$\Sp(8,\ZZ)$ whose fixed locus in $\HH_4$ is a divisor. This follows
easily from~\cite[Lemma 4.1]{Tai}: if an element
$\gamma\in\Sp(2g,\ZZ)$ of order~$m$ fixes $\tau\in\HH_g$ then it acts
on the tangent space with eigenvalues $e^{2\pi i(t_j+t_k)/m}$, where
$t_j$, $t_k\in\ZZ$ and $1\le j\le k\le g$. If $\tau$ is a general
point of a fixed divisor then $t_j+t_k\equiv 0\mod m$ for all but one
pair of indices, say $(j_0,k_0)$. But this is impossible if $g\ge
3$. To see this, we consider first the case $j_0=k_0$. We may assume
$j_0=1$, so $2t_1\not\equiv 0$, but then $t_1\equiv -t_2$ so
$2t_2\not\equiv 0$, and $(j_0,k_0)$ is not unique. On the other hand,
if $j_0\neq k_0$, we may assume $j_0=1$ and $k_0=2$, so
$t_1+t_2\not\equiv 0$; but in that case $t_3\equiv -t_2\equiv -t_1$ so
$2t_3\not\equiv 0$ and again $(j_0,k_0)$ is not unique.
\end{Proof}

\begin{cor}\label{Knef}
If $n\ge 3$, then the canonical bundle of $\AIGU4(n)$ is
nef. On the other hand, the canonical bundle of $\AVOR4(n)$ is never nef.
\end{cor}
\begin{remark}\label{models}
$\AIGU4(n)$ is a minimal model as defined in (\cite[Definition 2.13]{KM}),
because the singularities are terminal; but they are not $\QQ\,$-factorial
because $\Pi_2(4)$ is not simplicial, and some authors prefer to reserve
the term ``minimal model'' for the Mori category, whose objects are
projective varieties with $\QQ\,$-factorial terminal singularities. By toric
methods, following the argument of Fujino~\cite{Fuj}, a small
$\QQ\,$-factorialisation may be constructed, and this will be a
$\QQ\,$-factorial minimal model.
\end{remark}
\begin{theorem}\label{Kample}
If $n\ge 3$ then the canonical bundle of $\AIGU4(n)$ is ample.
\end{theorem}
\begin{Proof}
By Lemma~\ref{Katleveln} the canonical bundle satisfies the conditions
of Corollary~\ref{nefatleveln}, but with strict inequalities,
$a>12b/n>0$. Hence $K_{\AIGU4(n)}$ belongs to the interior of
$\Nef\big(\AIGU4(n)\big)\cap \big(\RR\,L+\RR\,\DIGU4(n)\big)$. Let
$H_0$ be an ample class on $\AIGU4(n)$ spanned by $L$ and $\DIGU4(n)$:
such an $H_0$ exists because $\AIGU4(n)$ is projective and
$\RR\,L+\RR\,\DIGU4(n)$ is the $\Sp(8,\ZZ/n)$-invariant part of $\Pic
\big(\AIGU4(n)\big)\tens\RR$, so if $H$ is some ample line bundle
class it is sufficient to take
$H_0=\sum_{\gamma\in\Sp(8,\ZZ/n)}\gamma(H)$. Now we copy the proof of
Kleiman's criterion given in~\cite[1.39]{KM}: $tK_{\AIGU4(n)}-H_0$ is
nef for $t\gg 0$, so for any dimension~$d$ subscheme
$Z\subset\AIGU4(n)$ we have $(tK_{\AIGU4(n)})^d\cdot Z\ge H_0^d\cdot
Z>0$ (this is a non-trivial step in the proof
of~\cite[1.38]{KM}). Therefore $tK_{\AIGU4(n)}$ is ample by the
Nakai-Moishezon criterion,~\cite[Theorem 1.37]{KM}.
\end{Proof}

Thus $\AIGU4(n)$ is the canonical model if~$n\ge 3$.

\section{The nef cone of the partial compactification}\label{nefconeofpartial}
We shall work with the partial compactification
$\cA'_g=\cA_g\cup D'_g$, sometimes with an additional level-$n$ structure
$\cA'_g(n)=\cA_g(n)\cup D'_g(n)$. If $n\ge 2$, then
$D'_g(n)=\sum_iD'_{g,i}(n)$ consists of several disjoint components, each of
which has a natural fibration $D'_{g,i}(n)\to \cA_{g-1}(n)$. For
$n\ge 3$ this is the universal family over $\cA_{g-1}(n)$, and for $n=1,2$
it is a family of Kummer varieties. Indeed
$$
D'_{g,i}(n)=\left(\ZZ^{2g-2}\rtimes \Gamma_{g-1}(n)\right)\backslash
\CC^{g-1}\times \HH_{g-1}.
$$
To describe this action let $m=(m', m'')$ with $m', m''\in \ZZ^{g-1}$
and $\gamma=\begin{pmatrix}A & B\\ C & D\end{pmatrix}\in
\Gamma_{g-1}(n)$. Then
$$
(m,\gamma)\colon (z,\tau)\mapsto\left((z+nm'+nm''\tau)(C\tau+D)^{-1},
(A\tau+B)(C\tau+D)^{-1}\right).
$$
If $n\ge 3$, then the fibre of the map $D'_{g,i}(n)\to \cA_{g-1}(n)$
over a point $[\tau]$ is the abelian variety $A_{n,n\tau}$ whose
period matrix is given by $(n\Eins_{g-1}, n\tau)$. For $n=1,2$ we
obtain the Kummer variety $A_{n,n\tau}/(\pm1)$.

Let $\Theta_{00}(z,\tau)\colon\CC^{g-1}\times\HH_{g-1}\to\CC$ be the
standard theta function. The automorphy factors of $\Theta^2_{00}$
define a $\QQ\,$-line bundle on $D'_{g,i}(n)$ which we shall denote by
$M'(n)$. For $n\ge 3$, let $N'=N_{D'_{g,i}(n)/\cA'_g(n)}$ be the
normal bundle of the boundary component $D'_{g,i}(n)$ in $\cA'_g(n)$.
\begin{lemma}\label{M'(n)}
If $n\ge 3$ then $ M'(n)=-nN'+L.$
\end{lemma}
\begin{Proof}
This is proved in \cite[Proposition 2.3]{Hu}. The proof consists of
comparing the cocycles of $M'(n)$ and $N'$.
\end{Proof}
\begin{proposition}\label{nefconeofA'}
The nef cone of $\cA'_g$ for $g\ge 2$ is given by
$$
\Nef(\cA'_g)=\left\{aL-bD'_g\mid b\ge 0, a\ge 12b\right\}.
$$
\end{proposition}
\begin{Proof}
The condition $b\ge 0$ is necessary, since $L$ is trivial on the
fibres of $D'_{g,i}(n)\to \cA_{g-1}(n)$, whereas $-D'_g(n)$ is ample
on the fibres (cf.~\cite[Proposition 1.8]{Mu}). In order to prove that
$a\ge 12b$ is a necessary condition we consider curves $C$ of the form
$X(1)\times\{A\}$ in $\cA'_g$ where $A$ is a fixed $(g-1)$-dimensional
abelian variety and $X(1)$ is the modular curve of level 1, i.e.\ we
consider a family of abelian varieties of type $E_{\tau}\times A$
where $E_{\tau}$ is an elliptic curve degenerating to a nodal
curve. Such a family is indeed contained in $\cA'_g$ and for
general~$A$, and $C.D'_g=1$. This is because the corresponding family
with a level-$n$ structure ($n\ge 3$ as usual) meets the boundary
transversally in a smooth point. Since the degree of the $\QQ\,$-line
bundle $L$ on $X(1)$ is $1/12$ we find the necessary condition $a\ge
12b$.

Next we shall prove that these conditions are sufficient. Here we
shall distinguish between curves $C$ which meet $\cA_g$ and curves $C$
which are contained in the boundary $D'_g$. For curves of the first
type the result was already proved in \cite[Proposition
1.4]{Hu}. Since the argument is very short we shall repeat it
here. Assume that $b\ge 0$ and $a\ge 12b$. Since $L$ is ample on the
Satake compactification, it follows that $L.C>0$, and we can assume
that $b>0$. Choose some $\varepsilon>0$ with $a/b > 12+\varepsilon$
and some point $[\tau]\in\cA_g$ on $C $. By a result of Weissauer
\cite[p. 220]{Wei} there exists a modular form $F$ of weight~$k$ and
vanishing order~$m$ such that $F(\tau)\neq 0$ and $m/k\ge
1/(12+\varepsilon)$. In terms of divisors this gives
$$
kL=mD'_g+D_F, \ C\not\subset D_F
$$
where $D_F$ is the closure in $\cA'_g$ of the divisor $\{F=0\}\subset 
\cA_g$. Hence
$$
\left(\sfrac k m L-D'_g\right).C=\sfrac 1 m D_F.C\ge 0
$$
and since $a/b>12+\varepsilon \ge k/m$ and $L.C>0$ we conclude that
$$
\left(\sfrac ab L-D'_g\right).C>\left(\sfrac km L-D'_g\right).C\ge 0.
$$
Finally let $C$ be a curve contained in $D'_g$. Here it is slightly easier to
work with level structures: we choose some $n\in\NN$ and assume that
$C\subset D'_{g,i}(n)$ for
some boundary component $D'_{g,i}(n)$ of $\cA'_g(n)$. By Lemma~\ref{M'(n)}
$$
\left(aL-bD'_{g,i}(n)\right)|_{D'_{g,i(n)}} =
\left(a-\sfrac bn\right) L+\sfrac bn M'(n).
$$
The condition $a\ge 12b$ for level~$1$ now becomes $a\ge 12b/n$. In
any case $a-b/n\ge 0$ and hence it suffices to prove that $M'(n).C\ge
0$. Fix a prime~$p$, and choose $n$ so that
$n\equiv 0\mod 4p^2$. If $m', m''\in \frac{1}{2p}\ZZ^{g-1}$
then the functions $\Theta^2_{m',m''}(z,\tau)$ define sections of $M'(n)$,
by \cite[Proposition 2.3]{Hu}: the proof uses the theta
transformation formula and the formulae ($\Theta$1)--($\Theta$3) from~\cite{Ig}
to show that $\Theta^2_{m',m''}(z,\tau)$ has the appropriate automorphy
factor. But this shows that $M'(n)$ is
generated by global sections and hence $M'(n).C\ge0$.
\end{Proof}

\section{Structure of the Voronoi boundary}\label{Vorboundary}
In this section we revert to the case $g=4$ and examine the geometry
of the Voronoi boundary in detail. Our chief purpose is to prove that
the fibration $D'_{4,i}(n)\to\cA_3(n)$ extends to the closure of
$D'_{4,i}(n)$ in the Voronoi compactification of $\cA_4(n)$. This
results in a fibration of each non-exceptional boundary divisor in
$\AVOR4(n)$ over $\AVOR3(n)=\AIGU3(n)$. The proof involves careful
study of the combinatorics of the cone $\Pi_2(4)$, and we also
assemble in this section some other results of that nature which we
shall need later.

\begin{proposition}\label{morphismexists}
Let $n \ge 3$, and suppose $\DVOR{4,i}(n)\subset\AVOR4(n)$ is the
closure of a boundary divisor, not contracted by
$\pi\colon\AVOR4(n)\to\AIGU4(n)$. Then there is a morphism
$$
p_i=p_{i,n}\colon \DVOR{4,i}(n)\To \AVOR3(n)
$$
extending the fibration $D'_{4,i}(n)\to\cA_3(n)$.
\end{proposition}

\begin{Proof}
We work, without loss of generality, with $\DVOR{4,1}(n)$,
corresponding to $\tau_{11}\to i\infty$: if $n\ge 3$ then $\DVOR{4,1}(n)$
is normal (see Remark~\ref{normalcpts}, below).
Thus we fix a rank~$3$ sublattice
$\LL_3=\ZZ x_2+\ZZ x_3+\ZZ x_4\subset\LL_4$ and set
$\MM_3=\Sym_2(\LL_3)$. The projection $\pr_1\colon\LL_4\to\LL_3$ with
kernel $\ZZ x_1$ induces a projection $\Sym_2\pr_1\colon\MM_4\to\MM_3$ with
kernel spanned by the $U^*_{1j}$, $1\le j\le 4$.

Consider the matrix
$$
\tilde\tau=\left(
\begin{array}{c|ccc}
\ast      & \tau_{12} & \tau_{13} & \tau_{14}\\
\hline
\tau_{12} & \tau_{22} & \tau_{23} & \tau_{24}\\
\tau_{13} & \tau_{23} & \tau_{33} & \tau_{34}\\
\tau_{14} & \tau_{24} & \tau_{34} & \tau_{44}
\end{array}
\right).
$$
Then
$$
z=(\tau_{12}, \tau_{13}, \tau_{14})\in\CC^3, \ \tau=(\tau_{ij})_{2\le
i,j\le 4}\in \HH_3.
$$
and the map $\tilde\tau\mapsto \tau$ is $\Sp(6,\ZZ)$-equivariant and
therefore induces a rational map $p_{1,n}\colon\DVOR{4,1}(n)\ratmap
\AVOR3(n)$. The problem is to extend this map to the cusps of
$\AVOR3(n)$.

We first check that $p_{1,n}$ extends over the smallest cusps, i.e.\ over
$\phi_n^{-1}(\cA_0)$, where $\phi_n\colon\AVOR4(n)\to\ASAT4$. This is the
only case which is nontrivial.
Near a component of $\phi_n^{-1}(\cA_0)$, the boundary divisor
$\DVOR{4,1}(n)$ is given by the fan $\Star\big(\Span{x_1^2},\Vor(4)\big)$
with respect to the lattice $\Bar\MM_4=\MM_4/\ZZ x_1^2$; see for
example~\cite[3.1]{Ful2}.
The map we are trying to extend, $p_{1,n}$, is given on the torus part of this
toric variety by forgetting all coordinates involving~$x_1$. More precisely,
$\DVOR{4,1}(n)$ is locally isomorphic to an analytic open set in
$$
X_{\Star\Span{x_1^2}}=
T_{\Bar\MM_4}\Temb\left(\Star\left(\Span{x_1^2},\Vor(4)\right)\right).
$$
The natural embedding $\MM_3\to\Bar\MM_4$ induces a map on the
corresponding tori
$$
T_{\Bar\MM_4}=\Hom(\Bar\MM_4, \CC^*)\To T_{\MM_3}
$$
which is $p_{1,n}$ on the torus part of $X_{\Star\Span{x_1^2}}$.

Now the result follows from Lemma~\ref{extendtoE} below. The extension
to the lower cusps, and the compatibility of the extensions, are
immediate consequences of the straightforward fact that if $\sigma$ is
an Igusa (i.e.\ Voronoi) cone in $\Sym_2^+(\ZZ^g)$ for $g<4$ then
$\Sym_2\pr_1(\sigma)$ is an Igusa cone in $\Sym_2^+(\ZZ^{g-1})$.
\end{Proof}
\begin{lemma}\label{extendtoE}
The map $T_{\Bar\MM_4}\to T_{\MM_3}$ extends to a
$\GL(\LL_3)$-equivariant map
$$
p_{1,n}\colon X_{\Star\Span{x_1^2}}\To T_{\MM_3}\Temb\big(\Vor(3)\big)
$$
of the corresponding torus embeddings.
\end{lemma}
\begin{Proof}
We need to check that the dual of the embedding, which may be thought
of as a projection $\Bar\MM_4\to\MM_3$ with kernel spanned by the
classes of $U^*_{1j}$, is a map of fans (the $\GL(\LL_3)$-equivariance
is automatic). To do that we must show that the projection of any
cone in $\Star\left(\Span{x_1^2}, \Vor(4)\right)$ lies in a cone of
$\Vor(3)$. By the definition of $\Star$, it is enough to show that if
$\sigma\in\Vor(4)$ and $\sigma\succ\Span{x_1^2}$, then $\Sym_2
\pr_1(\sigma)\subseteq \sigma'$ for some $\sigma'\in
\Vor(3)$. Moreover, since $\Vor(4)$ and $\Vor(3)$ are fans and
$\Sym_2\pr_1$ preserves the relation $\succ$ among cones, it is only
necessary to check this for top-dimensional cones in $\Vor(4)$ which have
$\Span{x_1^2}$ as a face. The result therefore follows from
Prop~\ref{exccones} and Prop~\ref{nonexccones}, below.
\end{Proof}

In verifying the assertion made in the above proof there are two cases
to be considered separately. If $\eta\prec \sigma$ (up to
$\GL(\LL_4)$-equivalence) then $\sigma$ corresponds to a point of the
exceptional locus $E\subset\AVOR4$. Otherwise $\sigma$ corresponds to
a point of $\AIGU4$.

\begin{proposition}\label{nonexccones}
Suppose that $\sigma\in\Vor(4)$ is of maximal dimension
(i.e.\ dimension~$10$), that $\Span{x_1^2}\prec\sigma$ and that no
$\GL(\LL_4)$-translate of $\eta$ is a face of~$\sigma$. Then
$\Sym_2\pr_1(\sigma)\in\Vor(3)$.
\end{proposition}

\begin{Proof}
In this case, $\sigma$ is equivalent under $\GL(\LL_4)$ to the first perfect
domain~$\Pi_1(4)$. (The level structure plays no role here.) More than that:
the subgroup of $\GL(\LL_4)$ that preserves $\Pi_1(4)$ permutes the
generating rays transitively, so $\sigma$ is even equivalent to $\Pi_1(4)$
under the stabiliser of $x_1^2$. To see that the rays are permuted
transitively, note first that the permutation matrices are in the stabiliser
of $\Pi_1(4)$ in $\GL(\LL_4)$, so all four monomial generators $x_i^2$ are
equivalent to one another and so are all six binomial generators
$(x_i-x_j)^2$. The element of $\GL(\LL_4)$ given by $x_i\mapsto x_i-x_2$
for $i\neq 2$ and $x_2\mapsto -x_2$ preserves $\Pi_1(4)$ but does not
preserve the distinction between monomial and binomial generators, so all the
generators are in one orbit.

Since, for any~$g$,
$$
\Pi_1(g)=\Span{x_1^2,\ldots,x_g^2, (x_i-x_j)^2\ (1\le i<j\le g)},
$$
the projection of $\Pi_1(g)$ to $\MM_{g-1}$ is $\Pi_1(g-1)$. Since
$\Pi_1(g)\in\Vor(g)$ for all~$g$ and $\Vor(g)$ is $\GL(\LL_g)$-invariant,
we certainly have $\Sym_2\pr_1\big(\Pi_1(g)\big)\in\Vor(g-1)$.
\end{Proof}

This part of the argument is not restricted to~$g=4$, but it only applies to
$\Pi_1(g)$. We want to mention an alternative proof, which uses the
information we have in a slightly different way. 

\begin{lemma}\label{dicinglemma}
Let $\LL$ be a lattice and $l_i\colon\LL\to\ZZ$ be linear forms such that the
quadratic form $\sum l_i^2$ is positive definite. Then the Delaunay
decomposition for the quadratic form $\sum \alpha_i l_i^2$ is independent
of the choice of positive constants $\alpha_i$ if, and only if, the forms
$l_i$ define a dicing; that is, the $0$-skeleton of the cell decomposition
defined by the hyperplanes $\{l_i(x)=n\}$ for $n\in\ZZ$ coincides with the
original lattice~$\LL$.
\end{lemma}
\begin{Proof}
\cite[Lemma 3.1]{ABH}.
\end{Proof}

Recall that it would be enough for our purposes to prove that
$\Sym_2\pr_1(\sigma)$ is contained in a cone of~$\Vor(3)$.
Every ray in $\Vor(4)$ is either a $\GL(\LL_4)$-translate of $\eta$ or
spanned by the square of a linear form. So if $\sigma\in\Vor(4)$ satisfies the
conditions of Proposition~\ref{nonexccones}, then $\sigma=\Span{\{l_i^2\}}$
for some linear forms $l_i\colon\LL_4^*\to\ZZ$. Now we can apply the following
proposition.

\begin{proposition}\label{dicing}
Suppose $\sigma\in\Vor(g)$ is a cone of maximal dimension which is spanned
by squares of linear forms. Then $\Sym_2\pr_1(\sigma)$ is contained in a cone
of~$\Vor(g-1)$.
\end{proposition}

\begin{Proof}
If $\sigma=\Span{\{l_i^2\}}$ then, since $\sigma$ is of maximal dimension,
$\sum l_i^2$ is positive definite.
Therefore, by Lemma~\ref{dicinglemma}
the $l_i$ define a dicing of $\LL_g^*\tens\RR$.
If $\xi'\in\LL_{g-1}\tens\RR$ is a
point of the $0$-skeleton of the decomposition induced by
the $\pr_1(l_i)$, then it is the projection of a cell in the dicing
of $\LL_g\tens\RR$ induced by the~$l_i$. Any vertex $\xi$ of this cell is
in~$\LL_g$, so $\xi'=\pr_1(\xi)$ is in~$\LL_{g-1}$. Therefore the
$\pr_1(l_i)$ induce a dicing of $\LL_{g-1}\tens\RR$.

The projection of any positive definite form is again positive
definite, so $\sum\left(\pr_1(l_i)\right)^2=\Sym_2\pr_1(\sum l_i^2)$
is positive definite. Therefore, again by Lemma~\ref{dicinglemma},
the Delaunay decompositions induced by any two forms in the interior
of $\Sym_2\pr_1(\sigma)$ are the same. Hence $\Sym_2\pr_1(\sigma)$ is
contained in a cone of~$\Vor(g-1)$.
\end{Proof}

Now suppose that $\eta\prec\sigma$, and that
$\Span{x_1^2}\prec\sigma$, so that $\sigma$ gives rise to a cone in
$\Star\big(\Span{x_1^2}, \Vor(4)\big)$. We need only consider
$10$-dimensional cones up to the action of the stabiliser $\tilde G_1$
in $\GL(\LL_4)$ of $\Span{x_1^2}$. Such a cone is spanned by $\eta$
and a $9$-dimensional facet of the second perfect
domain~$\Pi_2(4)$. These facets are described in~\cite{ER2}. The authors
of~\cite{ER2} have kept the coordinates $x_i$ and work with the cone
$\Psi'\big(\Pi_2(4)\big)$, but we prefer to work directly with $\Pi_2(4)$ and
to display the symmetry instead by using the coordinates $y_i=\Psi^{-1}(x_i)$
as in equation~\eqref{betarho} above. Facets of $\Pi_2(4)$ are then given by
setting some of the $\beta_{ij}$ and $\rho_{ij}$ equal to zero.

\begin{proposition}\label{orbitsofexccones}
Every $10$-dimensional cone $\sigma\in\Vor(4)$ with
$\Span{x_1^2,e}\prec\sigma$ is equivalent under
$\tilde G_1$ to one of the following three cones:
\begin{eqnarray*}
\Pi_2^1(4)&=&\{\beta_{14}=\beta_{34}=\rho_{13}=0\}+\eta;\\
\Pi_2^2(4)&=&\{\beta_{13}=\beta_{14}=\beta_{34}=0\}+\eta;\\
\Pi_2^3(4)&=&\{\beta_{14}=\beta_{34}=\rho_{24}=0\}+\eta.
\end{eqnarray*}
\end{proposition}
\begin{Proof}
Later (Corollary~\ref{exactorbits}) we shall show that $\Pi_2^1(4)$,
$\Pi_2^2(4)$ and $\Pi_2^3(4)$ are inequivalent under~$\tilde G_1$.
For now, since we are only interested in subcones of $\Pi_2(4)$, we need not
consider $\tilde G_1$ but only $G_1=\tilde G_1\cap G$, where
$G\subset \GL(\LL_4)$ is the subgroup that preserves~$\Pi_2(4)$.
Note that if $\sigma$ is as above, $g\in\tilde G_1$ and
$g(\sigma)\succ\Span{x_1^2,e}$ also, then $g\in G_1$ anyway. This is because
$g(e)$ is the barycentre of $g\big(\Pi_2(4)\big)$, so if $g$ does not
preserve $\Pi_2(4)$ then $e$ and $g(e)$ are in the interiors of different
top-dimensional cones of $\Igu(4)$ and cannot both be generators of
$g(\sigma)$, since $\Vor(4)$ is a refinement of $\Igu(4)$.

The symmetry group of $\Pi_2(4)$ is described
in~\cite{ER2}. It is a reflection group of order~$1152$, isomorphic to
the reflection group~$F_4$, generated by elements $k_i$, $(1\le i\le
4)$; $s_{ij}$, $(1\le i < j \le 4)$; and an extra transformation~$w$.
These are given by
$$
k_i(y_i)=-y_i,\quad k_i(y_j)=y_j\ (j\neq i);
$$
$$
s_{ij}(y_i)=y_j,\quad s_{ij}(y_j)=y_i,\quad s_{ij}(y_k)=y_k\ (k\neq i,j);
$$
$$
w(y_i)=-y_i+\half\sum_{k=1}^4 y_k.
$$
We claim that this group is~$G$; to show this, we must prove that it
is a subgroup of $\GL(\LL_4)$. Thus we need to check that the matrices
$\Psi^{-1}K_i\Psi$, $\Psi^{-1}S_{ij}\Psi$ and $\Psi^{-1}W\Psi$ are all
integral, where $K_i$, $S_{ij}$ and $W$ are the matrices of the above
transformations and $\Psi$ is the matrix of the Voronoi transformation
defined by equation~\eqref{Voronoitransformation}. Then
\begin{eqnarray*}
&&\Psi=\begin{pmatrix}
1 & 1 & 1 & 1\\
1 &-1 & 0 & 0\\
0 & 0 &-1 & 0\\
0 & 0 & 0 &-1
\end{pmatrix},\hbox{ so }
2\Psi^{-1}=\begin{pmatrix}
1 & 1 & 1 & 1\\
1 &-1 & 1 & 1\\
0 & 0 &-2 & 0\\
0 & 0 & 0 &-2
\end{pmatrix}\hbox{; and }\\
&&2W=\begin{pmatrix}
-1 & 1 & 1 & 1\\
 1 &-1 & 1 & 1\\
 1 & 1 &-1 & 1\\
 1 & 1 & 1 &-1
\end{pmatrix}
\end{eqnarray*}
so it is sufficient to check that $2\Psi^{-1}K_i\Psi$ and
$2\Psi^{-1}S_{ij}\Psi$ are congruent to zero mod~$2$ and that
$2\Psi^{-1}2W\Psi$ is congruent to zero mod~$4$. The first of these is
trivial since $K_i\equiv \Eins_4$ mod~$2$. For $S_{ij}$ it
is enough to notice that any two columns of $2\Psi^{-1}$ are equivalent
mod~$2$, so $2\Psi^{-1}S_{ij}\equiv 2\Psi^{-1}$ mod~$2$ and hence
$2\Psi^{-1}S_{ij}\Psi\equiv 2\Eins_4\equiv 0$. The case of~$W$ is
checked directly.

All these elements of~$G$ preserve $\eta$: they must do, as it is
spanned by the barycentre of~$\Pi_2(4)$. There are~$12$ rays
generating~$\Pi_2(4)$, spanned by $(y_i\pm y_j)^2$, and $G$ permutes
them transitively because $k_j\colon(y_i+y_j)^2\mapsto (y_i-y_j)^2$ and
$s_{ii'}s_{jj'}\colon(y_i+y_j)^2\mapsto(y_{i'}+y_{j'})^2$. Hence $G_1$, which
is the stabiliser of one of the rays (generated by
$x_1^2=(y_1+y_2)^2$) has order~$96$. The transformations $k_3$, $k_4$,
$k_1k_2$, $s_{12}$, $s_{34}$ and $w'=s_{14}s_{23}w$ all belong to~$G_1$;
and in fact they generate a group of order~$96$, which is therefore the
whole of~$G_1$. To see this, note that the elements $k_3$, $k_4$,
$k_1k_2$, $s_{12}$ and $s_{34}$ generate a group of
order~$32$, and the element $k_3w'$ has order~$3$: hence the group generated
has order at least~$96$.

The $G$-orbits of $9$-dimensional facets of $\Pi_2(4)$ are also
studied in~\cite{ER2}. There are exactly two such orbits, denoted RT
and BF. A facet is RT if it is $G$-equivalent to the facet
$\rho_{12}=\rho_{23}=\rho_{13}=0$, and BF if it is $G$-equivalent to
$\beta_{12}=\beta_{13}=\beta_{14}=0$: there are $16$~RT and $48$~BF
facets. The names come from the following representation, which will
also be useful to us. We construct a bicoloured graph on four vertices
numbered $1$~to~$4$: conventionally we think of these vertices as the
four corners of a square, numbered clockwise starting from the top
left. We join $i$ and $j$ with a red edge to represent the equation
$\rho_{ij}=0$ or with a black edge to represent $\beta_{ij}=0$. The
facets are then given by graphs with three edges that are forked
(there is a vertex of valency~$3$) or triangular. An RT facet is
$G$-equivalent to a facet described by a red triangular graph, and a
BF facet is $G$-equivalent to a facet described by a black forked
graph. The effect of $k_i$ on the graphs is to change the colour of
all edges having~$i$ as a vertex. $s_{ij}$ is just the transposition
$(ij)$ on the vertices. $w'$ interchanges a left black edge with a
right black edge (i.e.\ $\beta_{14}$ and $\beta_{23}$) and leaves
other black edges alone: to red edges it does the opposite, leaving
the left and right edges alone but interchanging top and bottom and
the two diagonals.

We are interested in facets adjoining $\Span{x_1^2}$ up to
$G_1$-equivalence. The coefficient associated to $x_1^2$ is
$\beta_{12}$, so we have $\beta_{12}\neq 0$: in other words, we look
only at graphs that do not have a black edge joining vertices~$1$
and~$2$.

There are $48$ facets adjoining $\Span{x_1^2}$, of which $12$ are~RT
and $36$ are~BF: this follows because $\Pi_2(4)$ has twelve edges, all
equivalent under~$G$, and each facet adjoins nine of them. $G_1$
preserves the property of being RT or BF, because $G$
does. $\Pi_2^1(4)$ is~RT, $\Pi_2^2(4)$ and $\Pi_2^3(4)$ are~BF.

The rest of the proof consists of checking that every facet of
$\Pi_2(4)$ that adjoins $\Span{x_1^2}$ occurs in one of these three
orbits, which is straightforward, using the description of the effects
of the generators above. The details are shown in Figure~1. Red edges
are shown as dotted lines and the vertices are numbered according to
the above convention, clockwise starting from the top left. Only half
the facets are shown, the others being their reflections (left-right)
under $s_{12}s_{34}$.
\end{Proof}
\begin{figure}[h!tp]
  \begin{center}
\epsfbox{facets.ps} 
  \end{center}
\begin{center}
Figure 1: Orbits of facets of the second perfect domain
\end{center}
\end{figure}

\begin{proposition}\label{exccones}
Each of the three cones $\Pi_2^i(4)$ projects under $\Sym_2\pr_1$ to a cone
contained in a cone of $\Vor(3)$.
\end{proposition}
\begin{Proof}
We simply check this for each case. (The representatives $\Pi_2^i(4)$
have been chosen so as to keep this part of the calculation fairly simple.)
Note that the projection of~$e$ is given by $\Sym_2\pr_1(e)=\bar e$, where
\begin{equation}\label{proje}
{\bar e}=(x_2-x_3)^2+(x_2-x_4)^2+x_3^2+x_4^2.
\end{equation}
Using this we have
\begin{eqnarray*}
\Sym_2\pr_1\big(\Pi_2^1(4)\big)&=&\Sym_2\pr_1
\left\langle x_1^2,x_2^2,x_4^2,(x_1-x_3)^2,(x_1-x_4)^2,\right.\\
&&\hskip-2pt\left.(x_2-x_3)^2,(x_2-x_4)^2,(x_3-x_4)^2,(x_1+x_2-x_3)^2,
e\right\rangle\\
&=&\left\langle x_2^2,x_4^2, x_3^2, x_4^2,(x_2-x_3)^2,(x_2-x_4)^2,(x_3-x_4)^2,
\right.\\
&&\left. (x_2-x_3)^2, (x_2-x_3)^2+(x_2-x_4)^2+x_3^2+x_4^2\right\rangle\\
&=&\Span{x_2^2,x_3^2,x_4^2,(x_2-x_3)^2,(x_2-x_4)^2,(x_3-x_4)^2}\\
&=&\Pi_1(3)
\end{eqnarray*}
and $\Pi_1(3)\in\Vor(3)$.

The first BF case is given by
\begin{eqnarray*}
\Sym_2\pr_1\big(\Pi_2^2(4)\big)&=&\Sym_2\pr_1
\left\langle x_1^2,x_2^2,x_3^2,x_4^2,(x_1-x_3)^2,(x_1-x_4)^2,\right.\\
&&\hskip-2pt\left.(x_2-x_3)^2,(x_2-x_4)^2,(x_3-x_4)^2,
e\right\rangle\\
&=&\left\langle x_2^2,x_3^2,x_4^2, x_3^2, x_4^2,(x_2-x_3)^2,
(x_2-x_4)^2,(x_3-x_4)^2,
\right.\\
&&\left. (x_2-x_3)^2+(x_2-x_4)^2+x_3^2+x_4^2\right\rangle\\
&=&\Span{x_2^2,x_3^2,x_4^2,(x_2-x_3)^2,(x_2-x_4)^2,(x_3-x_4)^2}\\
&=&\Pi_1(3)
\end{eqnarray*}

For the second BF case, $\Pi_2^3(4)$, we have
$$
\Sym_2\pr_1\big(\Pi_2^3(4)\big)=
\Span{x_2^2,x_3^2,x_4^2,(x_2-x_3)^2, (x_3-x_4)^2,{\bar e}}
$$
which by~\eqref{proje} is strictly contained in $\Pi_1(3)$.
\end{Proof}

\begin{cor}\label{exactorbits}
The $\tilde G_1$ orbits of $\Pi_2^1(4)$, $\Pi_2^2(4)$ and $\Pi_2^3(4)$ are
distinct.
\end{cor}

\begin{Proof} Since $\Pi_2^1(4)$ is an RT facet it is in a different
$G$-orbit from the other two, by~\cite{ER2}. As we have just seen,
$\Pi_2^2(4)$ projects onto a maximal-dimensional cone of $\Vor(3)$ and
$\Pi_2^3(4)$ does not, so they are inequivalent under~$\tilde G_1$.
\end{Proof}

Now we want to investigate the pullbacks of line bundles on $\AVOR3$ under
the morphisms~$p_i\colon \DVOR{4,i}(n)\To \AVOR3(n)$. We define
\begin{equation}\label{firstrestriction}
E_s(n)|_i=\DVOR{4,i}(n)\cap E_s(n).
\end{equation}
This intersection is either empty or a divisor on $\DVOR{4,i}(n)$
which is contracted to a variety of codimension~$\ge 2$ under the map
$\DVOR{4,i}(n)\to \DIGU{4,i}(n)$.

\begin{remark}\label{normalcpts}
The varieties $E_s(n)$ do not depend on $n$ for $n\ge 3$: more
precisely, all the $E_s(n)$ have the same normalisation (up to
isomorphism), independently of $s$ or $n$; and if $n\ge 3$ they are
normal. The normalisation is $\bar\cO(\eta)\subset
X_{\Vor(4)}$. Since the edges of $\Pi_2(4)$ are all equivalent
under~$G$ any two non-empty varieties $E_s(n)|_i$ are mutually
isomorphic as well. If $n\ge 3$ no nontrivial cone of~$\Vor(4)$ has
nontrivial stabiliser, since the principal congruence subgroup of
level~$n$ in $\GL(8,\ZZ)$ is torsion-free. This also implies that the
boundary divisors $\DVOR{4,i}(n)$ are normal.
\end{remark}

We recall here some facts about the structure of toroidal compactifications.
Recall from~\cite{Nam1} that any toroidal compactification of $\cA_g$ is a
disjoint union of strata of the form
$$
Z_{h,{\bar{\sigma}}}(n)=\cP_{g-h}(n)\backslash \HH_h \times
\CC^{h(g-h)}\times \cO({\bar{\sigma}})
$$
where $\cP_{g-h}(n)$ is a group which acts properly discontinously,
$\bar{\sigma}$ is a cone in some copy of $\Sym_2(\RR^{g-h})$ containing
some positive definite form, and $\cO(\bar\sigma)$ is the corresponding torus
orbit.

\begin{remark}\label{strata}
Suppose $C$ is an irreducible curve in $\AVOR4$. Then let $\sigma$ be
a maximal cone in $\Vor(4)$ such that $C$ is contained in the image in
$\AVOR4$ of the closure of the torus orbit $\cO(\sigma)$ (such a
$\sigma$ is unique up to the action of $\GL(\LL_4)$). If we assume
that $C$ is not contained in the exceptional divisor $E$, then
$\sigma$ must be of the form $\Span{l^2_1,\ldots, l^2_k}$ where the
$l_i$ are linear forms on $\LL_4^*$, as in the proof of
Proposition~\ref{dicing}.
\end{remark}

The connection with the strata $Z_{h,{\bar{\sigma}}}(n)$ is the following. Let
$$
U=\bigcap_{q\in\sigma} \Ker q\subset \MM_4^*\tens\RR,
$$
and set $h=\dim_\RR U$ and $V=\MM_4^*\tens\RR/U\cong \RR^{4-h}$. Then every
form $q\in \sigma$ defines a form $\bar{q}$ on $V$ and this defines an
injective map $\sigma\rightarrow\bar{\sigma}\subset \Sym_2^+(V)$.
\begin{lemma}\label{posforms}
$\bar{\sigma}\subset \Sym_2^+(V)$, as defined above, contains positive
definite forms.
\end{lemma}
\begin{Proof}
If $C\subset E$ then $e\in\sigma$ and since $e$~is positive definite
there is nothing to prove. Otherwise we prove this by induction on the
number~$m$ of generators of~$\sigma$. Suppose $C\not\subset E$ and
$\sigma=\Span{q_1,\ldots,q_m}$. We have $U=\bigcap\limits_{i=1}^m
\Ker q_i$, since if $q=\sum a_iq_i\in\sigma$ and $q_i(x)=0\in\MM_4\tens\RR$
for all~$i$ then $q(x)=0$. Thus $\bigcap \Ker \bar q_i=0$. Suppose
$0\neq x\in\Ker (\bar q_1+t\bar q_2)$ for some $t>0$. Then,
evaluating at $x\in V$, we get $\bar q_1(x)+t\bar q_2(x)=0$, and since both
forms are positive semidefinite this implies $\bar q_1(x) = \bar
q_2(x)=0\in\RR$. But since $\bar q_i$ is semidefinite, $\bar
q_i(x)=0$ if and only if $x\in \Ker \bar q_i$, so $x\in 
\Ker \bar q_1 \cap \Ker \bar q_2$. So $\Ker (\bar q_1+t\bar q_2)=\Ker
\bar q_1 \cap \Ker \bar q_2$, and this reduces to the case of
$\Span{q_1+t q_2,q_3,\ldots,q_m}$.
\end{Proof}

$C$ is then contained in a stratum of the form
$\cP_{4-h}(n)\backslash\HH_h\times\CC^{h(4-h)}\times
\cO(\bar{\sigma})$.

\begin{proposition}\label{pullbackD3}
Under the morphism $p_i$ of Proposition~\ref{morphismexists}
$$
p_i^*\big(D_3(n)\big)=\sum_{j\neq i}
\DVOR{4,j}(n)|_{\DVOR{4,i}(n)}+4\sum_s E_s(n)|_i.
$$
\end{proposition}
\begin{Proof}
The behaviour away from $\phi_n^{-1}(\cA_0)$ is clear and gives the
coefficient~$1$ for the boundary components~$\DVOR{4,j}$. It is
necessary to check the coefficient of~$E_s(n)$. The bundle
$p_i^*\big(D_3(n)\big)$ is given on $X_{\Star\Span{x_1^2}}$ by the support
function $\psi_3\circ\Sym_2\pr_1$, where $\psi_3$ is the support function on
$\Vor(3)$ that takes the value~$1$ on each primitive generator of a ray.
Hence the coefficient of~$E$ is $\psi_3(\bar e)=4$, by equation~\eqref{proje}.
\end{Proof}

We insert here some further details about the orbits of cones of $\Vor(4)$
that will be useful to us later on.

\begin{lemma}\label{dim2faces}
The dimension~$2$ faces of $\Pi_2(4)$ fall into two
orbits under the action of~$G$, the symmetry group of~$\Pi_2(4)$.
Representatives for these orbits are $\Span{x_1^2,x_2^2}$ and
$\Span{x_1^2,x_3^2}$.
\end{lemma}
\begin{Proof}
Any such face is equivalent under $G$ to a face spanned by $x_1^2$ and one
other generator of $\Pi_2(4)$. So we are interested in the $G_1$-orbits of
the other eleven generators. We can represent such a generator by a
bicoloured graph as we did for facets, only it is easier to use the
complementary graph, so that a red (respectively black) edge joining
vertices $i$ and $j$ represents $\rho_{ij}\neq 0$ (respectively
$\beta_{ij}\neq 0$). A generator of $\Pi_2(4)$ is thus represented by a
single edge. The generators of $G_1$ listed together with their action on
the graphs in the proof of Proposition~\ref{orbitsofexccones} all preserve
the property of an edge being horizontal. It is easy to see that the horizontal
and non-horizontal edges each form a $G_1$-orbit: see Figure~2.
The representatives given are defined by the non-vanishing of
$\rho_{12}$ and $\rho_{13}$ respectively.
\end{Proof}
\begin{figure}[h!tp]
  \begin{center}
\epsfbox{facets2.ps} 
 \end{center}
\begin{center}
Figure 2: Orbits of dimension~$2$ faces of the second perfect domain
\end{center}
\end{figure}

\begin{cor}\label{dim3cones}
Any cone of $\Vor(4)$ spanned by two rank~$1$ forms and a form of maximal
rank is $\GL(\LL_4)$ equivalent to one of the cones
$$
\sigma_3=\Span{x_1^2, x_2^2,e}\quad\hbox{  and  }\quad
\sigma'_3=\Span{x_1^2, x^2_2, e'}
$$
where
$$
e' =2(x_1^2+x_2^2+x_3^2+x_4^2+x_1x_3-x_1x_2-x_1x_4-x_2x_3-x_3x_4).
$$
\end{cor}
\begin{Proof}
Any such cone is equivalent to a cone spanned by~$e$ and a
dimension~$2$ face of $\Pi_2(4)$, so we can apply
Lemma~\ref{dim2faces}. However $\Span{x_1^2,x_2^2,e}=\sigma_3$, and
$\Span{x_1^2,x_3^2,e}$ is equivalent under $\GL(\LL_4)$
to~$\sigma'_3$; the element of $\GL(\LL_4)$ involved is simply the
transposition $x_2\leftrightarrow x_3$. Applying this to $e$ gives the
result.
\end{Proof}

\begin{lemma}\label{dim3faces}
The dimension~$3$ faces of $\Pi_2(4)$ fall into four orbits under the
action of~$G$, the symmetry group of~$\Pi_2(4)$. These orbits are to
be referred to as string, BF\/${}^*$, RT\/${}^*$ and disconnected: they
are represented by the cones $\Span{x_1^2,x_2^2, x_3^2}$,
$\Span{x_1^2,x_3^2, x_4^2}$, $\Span{x_1^2,x_4^2, (x_1-x_4)^2}$ and
$\Span{x_1^2,x_2^2, (x_3-x_4)^2}$ respectively.
\end{lemma}
\begin{Proof}
Such a face $\sigma$ is determined by three generators, i.e.\ by a bicoloured
graph with three edges. It is always possible to draw a forked or
triangular graph on the complement of such a graph. Therefore any
collection of three generators of $\Pi_2(4)$ spans a $3$-dimensional
face of $\Pi_2(4)$, since if we draw a triangular or forked graph on
the complement we specify, according to~\cite{ER2}, a facet containing
all those generators; and the facets, again according to~\cite{ER2},
are simplicial. (Note that the edges in the graphs in~\cite{ER2}
represent a condition $\beta_{ij}=0$ or $\rho_{ij}=0$, whereas for us
here they represent $\beta_{ij}\neq 0$ or $\rho_{ij}\neq 0$.)

If the graph representing $\sigma$ is itself triangular or forked,
then we may appeal directly to the argument of~\cite{ER2}. We conclude
that there are two orbits of these types, which we may call RT${}^*$ and
BF${}^*$, represented by the same graphs as the facets of types~RT
and~BF. Examples are $\Span{x_1^2,x_4^2, (x_1-x_4)^2}$ for~RT${}^*$
and $\Span{x_1^2,x_3^2, x_4^2}$ for~BF${}^*$. (A forked or triangular
graph is~RT if is triangular and has an even number of red sides:
otherwise it is~BF.)

So suppose that the graph is neither forked nor triangular (this means
that the rays not spanning~$\sigma$ do not span a facet
of~$\Pi_2(4)$). Then the vertices must have valencies $0$, $1$ or $2$,
and at least one of them has valency~$1$. The possibilities are that
all four vertices have valency~$1$; one has valency~$0$, one has
valency~$1$, and the other two have valency~$2$; or two have
valency~$1$ and two have valency~$2$. These possibilities are
illustrated in Figure~3.

We claim that the first two of these cases together form one
$G$-orbit, and that the third forms another. In the first case, of
valencies all equal to~$1$, the graph is a string. By applying
$s_{ij}$ we may assume that the string consists of edges joining $1$ to
$2$, $2$ to $3$ and $3$ to $4$. We may change the colour of the
outside edges by applying $k_1$ or $k_4$, and we may change the colour
of the central edge by applying $k_3k_4$, in each case without
changing anything else. Any graph of the second type (valencies $0$,
$1$, $2$, $2$) may be converted to a string by moving the double edge
to join $1$ to $2$ and then applying $w'$. So the first two types form
a single orbit. In the last type, where the graph is disconnected, we
may always move the double edge to join $1$ and $2$ and we may change
the colour of the remaining edge (necessarily joining $3$ and $4$) by
applying~$k_3$. It is also easy to see that this type cannot be
converted into a string.

Examples of these two possibilities, which we call ``string'' and
``disconnected'', are $\Span{x_1^2,x_2^2, x_3^2}$ and $\Span{x_1^2,x_2^2,
(x_3-x_4)^2}$ respectively.
\end{Proof}
\begin{figure}[h!tp]
  \begin{center}
\epsfbox{facets3.ps} 
 \end{center}
\begin{center}
Figure 3: Orbits of dimension~$3$ faces of the second perfect domain
\end{center}
\end{figure}
\begin{lemma}\label{mu6face}
Suppose $\sigma\prec\Pi_2(4)$ and that
$\sigma=\Span{l_1^2,l_2^2,l_3^2,l_4^2,l_5^2,l_6^2}$, and suppose that
the $l_i$ span a subspace of dimension~$3$. Then $\sigma$ is
$G$-equivalent to
$\Span{x_1^2,x_3^2,x_4^2,(x_1-x_3)^2,(x_1-x_4)^2,(x_3-x_4)^2}$.
\end{lemma}
\begin{Proof}
Without loss of generality we may assume that $l_1$, $l_2$ and $l_3$
are linearly independent. Since every face of $\Pi_2(4)$ is
simplicial, $\Span{l_1^2,l_2^2,l_3^2}\prec\Pi_2(4)$, so according to
Lemma~\ref{dim3faces} it is equivalent to either $\Span{x_1^2,x_2^2,
x_3^2}$ or $\Span{x_1^2,x_3^2, x_4^2}$. Type RT${}^*$ is excluded by
the linear independence condition, and disconnected type is excluded
because for any other generator $l_4^2$ of $\Pi_2(4)$ the linear forms
$x_1$, $x_2$, $x_3-x_4$ and $l_4$ span a space of dimension~$4$. We
prefer to replace $\Span{x_1^2,x_2^2, x_3^2}$ by the equivalent face
$\Span{x_1^2,x_3^2, (x_3-x_4)^2}$ (also of string type -- apply~$w'$
followed by $s_{12}s_{34}$).

Now the result follows from the observation that any seven linear
forms whose squares are generators of $\Pi_2(4)$ span a linear space
of dimension~$4$. Therefore the six $l_i$ are all the linear forms
whose squares are generators of $\Pi_2(4)$ and which lie in the linear
span of $l_1$, $l_2$ and $l_3$. In both cases this gives 
$\Span{x_1^2,x_3^2,x_4^2,(x_1-x_3)^2,(x_1-x_4)^2,(x_3-x_4)^2}$.
\begin{figure}[h!tp]
  \begin{center}
\epsfbox{facets4.ps} 
 \end{center}
\begin{center}
Figure 4: An orbit of faces of the second perfect domain spanned by
six generators
\end{center}
\end{figure}
\end{Proof}
\begin{remark}\label{opposites}
The graph corresponding to this example is a coloured complete graph,
as shown on the left of Figure~4. Applying $w'$ gives a graph which is
a bicoloured triangle. Notice that each edge of this graph has a
distinguished opposite edge, which shares an even number of vertices with
it but is of the opposite colour: this will be used below in
Proposition~\ref{mu6}.
\end{remark}

\section{Proof of the main result}\label{mainsection}
Recall that $\DVOR{g}(n)$ denotes the closure in $\AVOR{g}(n)$ of the
boundary $D'_g(n)$. We refer to the irreducible components of $\DVOR{g}(n)$
as $\DVOR{g,i}(n)$, or simply as $D_{g,i}(n)$ if $g\le 3$ (when there is
only one toroidal compactification we need consider). On $\DVOR{g,i}$ we
define the line bundle
\begin{equation}\label{defineM}
M_{g,i}(n)=-nN_{\DVOR{g,i}(n)/\AVOR{g}(n)}+L.
\end{equation}
Clearly, in view of Lemma~\ref{M'(n)}, $M_{g,i}(n)$ is an extension of the
line bundle $M'(n)$ introduced in section~II.

In the case $g=4$ the boundary of $\AVOR4(n)$ decomposes as
$$
\AVOR4(n)\setminus\cA_4(n)=\DVOR4(n)+E(n)=\sum_i\DVOR{4,i}(n)+\sum_s E_s(n).
$$
We define the line bundle on~$\DVOR{4,i}$
\begin{equation}\label{defineJ}
J_4(n)=J_{4,i}(n)=M_{4,i}(n)-\sum_s n E_s(n)|_i.
\end{equation}

The proof of Proposition~\ref{nefconeofA'} shows that certain theta
functions define sections of $M'_{4,i}(n)$. The first technical result
of this section, Proposition~\ref{sectionsofJ}, is that these
sections extend to sections of $J_4(n)$.

The following notation will be used throughout the rest of the paper: if $I$
is a set of indices then we write $\DVOR{g,I}(n)$ for
$\bigcap\limits_{i\in I}\DVOR{g,i}(n)$, and if $\cF$ is a bundle (or sheaf,
etc.) on some variety containing $\DVOR{g,I}(n)$ we denote the restriction
$\cF|_{\DVOR{g,I}(n)}$ more simply by $\cF|_I$. We have already used this
convention above (equation~\eqref{firstrestriction}): as we did there,
we normally
abuse notation by writing $\DVOR{4,ij}(n)$ and $\cF|_{ij}$ rather than
$\DVOR{4,\{i,j\}}(n)$ and $\cF|_{\{i,j\}}$, etc.

\begin{proposition}\label{sectionsofJ}
Let $p$ be a prime and $n\equiv 0$ mod~$4p^2$. If the characteristics
$m', m'',\bar{m}', \bar{m}'' \in
\frac {1}{2p} \ZZ^{g-1}$, then the functions
$\Theta_{m',m''}(z,\tau)\Theta_{\bar{m}',
\bar{m}''}(z,\tau)$ define sections of the line bundle $J_4(n)$.
\end{proposition}
\begin{Proof}
Since the group $\Sp(8,\ZZ/n)$ acts transitively on the boundary
components $\DVOR{4,i}(n)$ we can again restrict ourselves to the
standard boundary component $\DVOR{4,1}(n)$, given by the line
$\QQ\, \bde_1$ in $\QQ^8=\QQ\,\bde_1\oplus\cdots\oplus\QQ\,\bde_8$.
We write $\DVOR{4,1}(n)\cap\cA'_4(n)=D'_{4,1}(n)$. We have already observed
that the functions
$\Theta_{m'm''}(z,\tau)\Theta_{\bar{m}'\bar{m}''}(z,\tau)$ have the
correct transformation behaviour. We shall have to study how the
sections defined by these functions extend to the generic point of the
intersections $\DVOR{4,1j}(n)=\DVOR{4,1}(n)\cap \DVOR{4,j}(n)$, $j\neq 1$
and to the generic point of the divisors~$E_s(n)|_1$.

A standard calculation shows that $\Sp(8,\ZZ/n)$ acts transitively on pairs
$\big(\DVOR{4,i}(n), \DVOR{4,j}(n)\big)$ with
$\DVOR{4,ij}(n)\neq\emptyset$. Hence we can work
with the standard cusp corresponding to the isotropic subspace
$\QQ\,(\bde_1\wedge \bde_2\wedge \bde_3\wedge \bde_4)$ and we
can, moreover, take $j=2$ and assume that $\DVOR{4,1}(n)$ and $\DVOR{4,2}(n)$
correspond to the rays $\Span{x^2_1}$ and $\Span{x^2_2}$ in $\MM_4\tens\RR$.
The cone $\Span{x^2_1,x^2_2}\in\Vor(4)$ has dual cone given by
$$
\Span{x_1^2,x_2^2}^{\vee}=
\Span{U_{11}, U_{22}, \pm U_{33}, \pm U_{44}, \pm U_{ij}\ (i\neq j)},
$$
where $\{U_{ij}\}$ is the dual basis to $\{U_{ij}^\ast\}$.

We have the partial quotient
\begin{eqnarray*}
\HH_4    &\To &\CC\times\CC\times(\CC^*)^8=
T_{\MM_4}\Temb\left(\Span{x_1^2,x_2^2}\right)\\
(\tau_{ij}) &\longmapsto&(t_{11}=e^{2\pi i\tau_{11}/n},
t_{22}=e^{2\pi i\tau_{22}/n}, t_{ij}=e^{2\pi i\tau_{ij}/n}).
\end{eqnarray*}
in which $\DVOR{4,1}(n)$ corresponds to $\{t_{11}=0\}$ and $\DVOR{4,12}(n)$
corresponds to $\{t_{11}=t_{22}=0\}$.

We now have to study the theta functions $\Theta_{m' m''}(z,\tau)$
where $m', m''\in(1/2 p)\ZZ^{g-1}$. The transformation of these
functions with respect to $\Sp(8,\ZZ)$ is given by the theta
transformation formula \cite[pp. 84, 85]{Ig}. Note that
characteristics of the form $(m', m'')$ with $m',
m''\in\frac{1}{2p}\ZZ^{g-1}$ are transformed to characteristics of the
same type. 

The connection between the variables
$(z,\tau)\in\CC^3\times\HH_3$ and $\HH_4$ is the following.
Recall that $\DVOR{4,1}(n)$ corresponds to $\tau_{11}\to i\infty$.
In terms of the coordinates $t_{ij}=e^{2\pi i\tau_{ij}/n}$ we have
\begin{eqnarray}\label{thetaformula}
\lefteqn{\Theta_{m' m''}(z, \tau)=}\nonumber\\
&\ & \left(\Sum_{q\in\ZZ^3} t^{\half (q_2+m'_2)^2n}_{22}
t^{\half(q_3+m'_3)^2n}_{33} t^{\half(q_4+m'_4)^2n}_{44}\right)\cdot\\
&&\prod_{2\le i < j\le 4} t_{ij}^{(q_i+m'_i)(q_j+m'_j)n}
t_{12}^{(q_2+m'_2)n}
t_{13}^{(q_3+m'_3)n} t_{14}^{(q_4+m'_4)n} e^{2\pi i^t(q+m')m''}.\nonumber
\end{eqnarray}
The claim that the sections of $M'_{4,i}(n)$ defined by the
products of theta functions
$\Theta_{m' m''}(z,\tau)\Theta_{\bar{m}'\bar{m}''}(z,\tau)$ can be extended
over $\DVOR{4,12}(n)$ now follows from the
observation that the exponent of $t_{22}$ in $\Theta_{m'm''}(z,\tau)$ is
equal to $(q_2+m'_2)^2 n/2$ and, in particular, non-negative.

Next we study the extension to the generic point of a divisor
$E_s(n)|_1$. Again we claim that $\Sp(8,\ZZ/n)$ acts transitively on the
pairs $(\DVOR{4,i}(n), E_s(n)|_i)$. By the action
of the group $\GL(\LL_4)$ on $\MM_4$ we can assume that $E_s(n)$ corresponds
to the central ray~$\eta$ in the second perfect cone. The transitivity now
follows from the observation at the start of the proof of
Proposition~\ref{orbitsofexccones} that
$G$, the stabiliser of $\Pi_2(4)$ in $\GL(4,\ZZ)$, permutes the
generators of $\Pi_2(4)$ transitively and preserves~$\eta$.

We can, therefore, restrict our attention to
the $2$-dimensional cone $\Span{x^2_1, e}$. But
$$
e=2U_{11}^\ast+2U_{22}^\ast+2U_{33}^\ast+2U_{44}^\ast+U_{12}^\ast
-U_{13}^\ast-U_{14}^\ast-U_{23}^\ast-U_{24}^\ast
$$
and, using this, a straightforward calculation shows that the dual cone
is given by
\begin{eqnarray*}
\Span{x_1^2,e}^{\vee}&=&\left\langle(U_{11}-2U_{12}), U_{12},
\pm (U_{22}-U_{33}), \pm (U_{22}-U_{44})\right.\\
&&\left.\pm(U_{24}+U_{12}),\pm (U_{24}-U_{13}),\pm (U_{24}-U_{14}),
\pm (U_{24}-U_{23})\right.\\
&&\left.\pm (2U_{12}-U_{22}),\pm U_{34}\right\rangle.
\end{eqnarray*}
Hence $T_{\MM_4}\Temb\big(\Span{x_1^2,e}\big)\cong\CC^2\times (\CC^*)^8$
and the torus embedding is given by
\begin{eqnarray*}
T_{\MM_4}&\!\To\!&T_{\MM_4}\Temb\left(\Span{x_1^2,e}\right)
\cong\CC\times\CC\times(\CC^*)^8\\
(t_{ij})&\!\longmapsto\!& 
(t_{11} t^{-1}_{12}, t_{12}, t_{22} t^{-1}_{33}, t_{22}
t_{44}^{-1}, t_{24} t_{12}, t_{24} t^{-1}_{13},
t_{24} t^{-1}_{14}, t_{24} t^{-1}_{23}, t^2_{12} t^{-1}_{22}, t_{34}).
\end{eqnarray*}
Let $T_1,\ldots, T_{10}$ be the obvious coordinates on $\CC^2\times (\CC^*)^8$,
corresponding to $U_{11}-2U_{12}, U_{12}, U_{22}-U_{33},\ldots,
2U_{12}-U_{22},U_{34}$. The hyperplane $\{T_1=0\}$ describes
$\DVOR{4,1}(n)$ and $\{T_1=T_2=0\}$ defines $E_s(n)|_1$. A
straightforward calculation shows
\begin{eqnarray}\label{dualx1}
&&t_{11}= T_1T_2^2,\quad t_{12}=T_2, \quad t_{13}=T^{-1}_2 T_5 T^{-1}_6,
\quad t_{14}=T^{-1}_2 T_5 T^{-1}_7,\nonumber\\
&&t_{22}=T^2_2 T^{-1}_9,\quad t_{23}=T^{-1}_2 T_5 T^{-1}_8,\quad 
t_{24}=T^{-1}_2 T_5,\\
&&t_{33}=T^2_2 T^{-1}_3 T^{-1}_9,\quad t_{34}=T_{10},\quad 
t_{44}=T^2_2 T^{-1}_4 T^{-1}_9.\nonumber
\end{eqnarray}
Combining this with formula~\eqref{thetaformula}, we can write the
function $\Theta_{m',m''}(z,\tau)$ in terms of these new coordinates
$T_1,\ldots,T_{10}$. We are interested in the exponent of $T_2$. If
$$
w_2=(q_2+m'_2), \ w_3=(q_3+m'_3),\ w_4=(q_4+m'_4)
$$
then this exponent is given by
$$
\ell(w_2,w_3,w_4)=n(w_2-w_2w_4-w_3-w_4-w_2w_3+w_2^2+w_3^2+w_4^2).
$$
Another straightforward calculation shows that this function assumes
its minimum for $w_2=0$, $w_3=w_4=1/2$, where we find that $\ell(0, 1/2,
1/2)=-n/2$. Altogether for products of the form $\Theta_{m'
m''}(z,\tau)\Theta_{\bar{m}'\bar{m}''}(z,\tau)$ we pick up poles of
order at most $n$. On the other hand $t_{11}=T_1T^2_2$ shows that we
have a zero of order $2n$ and this means that, in total, we have a
zero of order at least $n$.
\end{Proof}

Before we give the proof of the main theorem we want to introduce the notion
of depth of an irreducible curve~$C$. Recall that we have a morphism
$$
\phi_n\colon \AVOR4(n)\To \AVOR4
\To \ASAT4=\cA_4 \amalg\cA_3\amalg\cA_2\amalg\cA_1\amalg \cA_0.
$$
The depth of an irreducible curve $C\subset\AVOR4(n)$ is
defined by
\begin{equation*}
\depth(C):=\min\{k\mid \phi_n(C)\cap\cA_{4-k}\neq\emptyset\}.
\end{equation*}
Obviously $0\le \depth(C)\le 4$ and $\depth(C)=0$ if and only if $C$
is not contained in the boundary. In the rest of the paper we shall treat
each case in turn, starting with depth~$4$
(subsection~\ref{depth4section}) and then going on to depth~$0$ in
subsection~\ref{depth0section}, depth~$1$ in
subsection~\ref{depth1section}, depth~$2$ in subsection~\ref{depth2section}
and finally depth~$3$ in subsection~\ref{depth3section}.

\subsection{Curves of depth 4}\label{depth4section}

For a depth~$4$ curve, the question of whether it meets a given divisor
negatively is a purely toric one, depending only on facts about~$\Vor(4)$.

\begin{proposition}\label{nefondepth4}
A divisor $aL-bD_4-cE$ on $\AVOR4$ has non-negative intersection
with all irreducible curves $C$ of depth~$4$ if and only if $b\ge 2c\ge 0$.
\end{proposition}

\begin{Proof}
Since $C$ is mapped to a point in the Satake compactification it
follows that $L.C=0$. First we suppose that $C=C_E$ is a curve in
$E$. It is known that $\AVOR4$ is projective: this was proved by
Alexeev in~\cite{Al} for all $\AVOR{g}$ but seems to have been known
for much longer for $g=4$; see for instance~\cite{Nam2}. Therefore
there is an ample line bundle on $E$ which is the restriction of a
line bundle on~$\AVOR4$. But we saw in Proposition~\ref{picard} that
$\Pic(\AVOR4)$ is generated by $L$, $\DVOR4$ and~$E$; and the first
two of these are pulled back from $\AIGU4$ and hence zero on~$E$. So
either $E|_E$ or $-E|_E$ is ample, and it is easy to see that in fact
$-E|_E$ is ample. This follows because $\Pi_2(4)$ is contained in a
half-space. Hence on the toric variety $T_{\MM_4}\Temb\big(\Vor(4)\cap
\Pi_2(4)\big)$ got by subdividing the second perfect domain into
Voronoi cones, with torus-invariant divisors $E=\bar\cO(\eta)$ (this
is an abuse of notation as it is the normalisation
of~$E\subset\AVOR4$) and $D_1,\ldots,D_{12}$ given by the generators,
there is a linear relation involving $E$ and all the $D_i$ with
positive coefficients. So, on the toric variety, $-E|_E$ is effective;
and this remains true on~$\AVOR4$. So
$H.C\ge 0$ if and only if~$c\ge 0$.

Actually we can do better. If we work instead with
$\Psi'\big(\Pi_2(4)\big)$ and use the linear relation induced by the
linear form~$\sum U_{ii}$, we see that on the toric variety
$-4E|_E=K|_E$, so that $E\subset X_{\Vor(4)}$ is a toric Fano variety.

Next we consider the case of a curve $C$ of depth~$4$ that does not
meet~$E$. Then $C$ is contained in the boundary $\DIGU4$ of $\AIGU4$,
and the same considerations as above, applied to $\AIGU4\to\ASAT4$,
show that $-\DIGU4|_{\DIGU4}$ is ample.

We remark that the morphism $\AVOR4\to\ASAT4$ is a normalised blow-up
of some sheaf of ideals, (see~\cite{Ch} and \cite[Section IV]{SC}),
and this is also sufficient for our purposes.

For other curves of depth~$4$, it is convenient to work on $\AVOR4(n)$
for some $n\ge 3$: pulling back by $\alpha_{n,\Vor}$, we must show
that $\big(bD_4(n)+cE(n)\big).C\le 0$ for every irreducible curve $C$
of depth~$4$ if and only if $b\ge 2c\ge 0$. Notice that this will also
prove that these conditions $b\ge 2c\ge 0$ are necessary for
$aL-bD_4-cE$ to be nef, as claimed in~Theorem~\ref{mainthm}.

It is enough to consider the curves $C$ corresponding to
codimension~$1$ cones $\sigma\in\Vor(4)$. This is because we need only
consider irreducible curves, and any such curve in $E(n)$ lifts to a
single irreducible component $E$ of the boundary of $X_{\Vor(4)}$. But
such a component is a toric variety, and the curves corresponding to
codimension~$1$ cones generate the cone of effective curves in
$E$. Indeed, such curves generate the whole of $A_1(E)$, by for
instance~\cite[Proposition 10.3]{Dan}, and rational equivalence
implies numerical equivalence.

We have already dealt with such curves in the case where
$\sigma\succ\Span{g(e)}$ for some $g\in\GL(\LL_4)$, because they are
contained in~$E(n)$. So it remains to deal with
$\sigma\prec\Pi_2(4)$. Up to $\GL(\LL_4)$-action there are two such
cones (RT and BF facets in the notation of~\cite{ER2}).
We choose to work with the cones
$\sigma_0=\{\beta_{13}=\beta_{14}=\rho_{34}=0\}$ and
$\sigma_1=\{\beta_{13}=\beta_{14}=\beta_{34}=0\}$. Each of these is a
$9$-dimensional face of the second perfect cone $\Pi_2(4)$ and defines
a rational curve $C\cong\PP^1$ in $X_{\Vor(4)}$. 

Let $\Pi_2^1(4)'$ be the $10$-dimensional cone $\Span{\sigma_0,e}$. In
the language of~\cite{ER2}, $\Pi_2^1(4)'$ is a type~III domain: in the
classification of Proposition~\ref{orbitsofexccones} it is equivalent
to $\Pi_2^1(4)$. The facet $\sigma_0$ is an RT facet
of~$\Pi_2(4)$. The transformation $x_1\leftrightarrow x_3,
x_2\leftrightarrow x_4$ leaves $\sigma_0$ invariant, but maps $\eta$
to a ray $\eta'=\Span{e'}$ and $\Pi_2^1(4)'$ to another
$10$-dimensional cone, a part of a translate of $\Pi_2(4)$, which is
again a type~III domain. 

The geometric situation is this: the curve $C_0\subset X_{\Vor(4)}$
corresponding to
$\sigma_0$ is contained in nine boundary components, which we call
$\DVOR{4,2},\ldots,\DVOR{4,10}$ (here $10$ means `ten', not $\{1,0\}$:
the reason for the indexing will appear
below), out of the twelve boundary components
$\DVOR{4,1},\ldots,\DVOR{4,12}$ corresponding to $1$-dimensional faces of
$\Pi_2(4)$. These are the ones belonging to the $1$-dimensional faces
of $\sigma_0$. It is met (transversely) by two exceptional divisors
$E$, $E'$ corresponding to the rays $\eta$ and $\eta'$. No other
invariant divisors meet~$C$, and because of the level structure $E$
and $E'$ give distinct disjoint components of~$E(n)$.

The form $\half\Tr'$ which we introduced in the proof of
Proposition~\ref{D4versusD4Vor} takes the value~$1$ on each primitive
generator of $\sigma_0$, while $\half\Tr'(e)=4$ and
$\half\Tr'(e')=5$. This shows that
$$
\DVOR{4,2}+\ldots+\DVOR{4,10}+4E+5E'+R\sim 0
$$
where $R$ is a divisor in $X_{\Vor(4)}$ which does not meet~$C_0$.
This implies
$$
D_4(n).C_0=(\DVOR{4,2}+\ldots+\DVOR{4,10}+4E+4E').C_0=-E'.C_0=-1.
$$
We also have $L.C_0=0$ and $E(n).C_0=(E+E').C_0=2$. Hence
$$
\big(-bD_4(n)-cE(n)\big).C_0=b-2c\ge 0,
$$
which is the desired inequality.

Finally we do the same calculation for~$\sigma_1$. This cone has
$\Span{\sigma_1,e}=\Pi_2^2(4)$ which is a type~II domain: $\sigma_1$
itself is a BF facet, since
$w(\sigma_1)=\{\beta_{12}=\beta_{23}=\beta_{24}=0\}$. It forms the
boundary between $\Pi_2(4)$ and $\Pi_1(4)$ and shares an
$8$-dimensional face with $\sigma_0$.

Choosing the numbering suitably, we have the following geometric
picture: the curve $C_1$ lies in the intersection of the nine boundary
divisors $\DVOR{4,1},\DVOR{4,2},\ldots,\DVOR{4,9}$, and is met transversely by
$E$ and another boundary component $\DVOR{4,0}$ corresponding to
$\Span{(x_1-x_2)^2}$. (With this choice of indexing the $1$-dimensional
faces of $\Pi_1(4)$ correspond to $\DVOR{4,0},\ldots,\DVOR{4,9}$.) Now
$\half \Tr'\big((x_1-x_2)^2\big)=2$ so
$$
\DVOR{4,1}+\cdots+\DVOR{4,9}+4E+2\DVOR{4,0}+R\sim 0
$$
for some $R$ not meeting $C_1$; and hence
$$
D_4(n).C_1=(\DVOR{4,0}+\DVOR{4,1}+\cdots+\DVOR{4,9}+4E).C_1=-\DVOR{4,0}.C_1=-1.
$$
We also have $E(n).C_1=E.C_1=1$ so 
$$
\big(-bD_4(n)-cE(n)\big).C_1=b-c\ge b-2c\ge 0,
$$
which completes the proof.
\end{Proof}

\subsection{Curves of depth 0}\label{depth0section}

The method we use in this case, of curves that are not contained in the
boundary of $\AVOR4$, is analogous to the proof of
Proposition~\ref{nefconeofA'}, but rather more complicated. We have to
produce a modular form vanishing to sufficiently high order along each
boundary component of $\AVOR4$. It is relatively easy to supply the
modular form: our proof that it does indeed have the required vanishing
is neither simple nor elegant.

\begin{proposition}\label{nefondepth0}
Let $C$ be a depth~$0$ curve and let $H=aL-bD_4-cE$ be a divisor
on $\AVOR4$ with $a\ge 0$, $a-12b\ge 0$ and $b\ge 2c\ge 0$. Then $H.C\ge 0$.
\end{proposition}

\begin{Proof}
It is simpler to work this time with $\DVOR4$ rather than $D_4$, so we
write $H=\alpha L-\beta\DVOR4-\gamma E$ (as in Remark~\ref{nefbyVor})
and assume that $\beta\ge 0$, $\alpha -12\beta\ge 0$ and $\gamma\ge
4\beta\ge \sfrac{8}{9}\gamma$.

We note that $L.C>0$ since $C$ maps to a curve in the Satake compactification
and $L$ is ample on $\ASAT4$. It is enough to prove that
$H.C>0$ if $a-12b>0$. Choose some $\varepsilon>0$ with
$a/b>12+\varepsilon$, and let $F$ be a modular form of weight $k$,
with $F|_C\not\equiv 0$, vanishing of order $m\ge k/(12+\varepsilon)$
on $\DVOR4$ and order $r\ge 9m/2$ on~$E$: such a form exists by
Proposition~\ref{9/2-vanishing}, below.

Now we can write
$$
kL=m\DVOR4+rE+D_F, \ C\not\subset D_F
$$
where $D_F$ is the zero divisor of $F$ on $\AVOR4$ (that is, the
closure in $\AVOR4$ of the set $\{F=0\}\subset\cA_4$). So 
$$
\left(\sfrac{k}{m} L-\DVOR4-\sfrac{r}{m}E\right).C=\sfrac{1}{m} D_F.C\ge 0.
$$
Since $a/b>12+\varepsilon\ge k/m$ and $r/m\ge 9/2$, and $E.C\ge 0$
since $C$ is not of depth~$4$, it follows that
$$
\left(\sfrac{\alpha}{\beta} L-\DVOR4-\sfrac{\gamma}{\beta}E\right).C
>\left(\sfrac{k}{m} L-\DVOR4-\sfrac{r}{m}E\right).C\ge 0.
$$
\end{Proof}
It remains to establish that the modular form~$F$ exists.

\begin{proposition}\label{9/2-vanishing}
Given an irreducible curve $C\subset\AVOR4$ of depth~$0$, There exists
a $k\in\NN$ and a modular form $F$ for $\Sp(8,\ZZ)$ of weight $k$ with
$F|_C\not\equiv 0$, such that $F$ vanishes of order~$m$ on $\DVOR4$
and order~$r$ on~$E$, and $m/k\ge 1/(12+\varepsilon)$ and $r\ge 9m/2$.
\end{proposition}

\begin{Proof}
Apart from the inequality $r\ge 9m/2$ this is the result of Weissauer
\cite[p. 220]{Wei} that we used in the proof of
Proposition~\ref{nefconeofA'}. We shall prove that the forms that
Weissauer constructs also fulfill the inequality $r\ge 9m/2$. 
For this purpose we need to recall his construction.

Let $l=2p$ with $p$ prime and consider the set $\cM$ of all
characteristics in $(\sfrac{1}{l}\ZZ/\ZZ)^{8}$ of the form
$$
m=(m_{(p)}, m_{(2)})\in 
(\sfrac{1}{p}\ZZ/\ZZ)^{8}\oplus(\hf \ZZ/\ZZ)^{8},
\ m_{(p)}\neq 0.
$$
For a characteristic $m=(m', m'')$ with $m', m''\in\RR^4$ the associated theta
constant is defined by
$$
\Theta_m(\tau, 0)=\sum_{q\in\ZZ^g} 
e^{2\pi i[\half(q+m')\tau^t(q+m')+(q+m')^tm'']}
$$
For $\tilde\cM\subset \cM$, define
$$
\Theta_{\cM,\tilde\cM} (\tau)=\prod_{m\in\cM\setminus \tilde\cM}
\Theta_m(\tau, 0)^l
$$
and
$$
F_r(\tau)=\sum_{M\in\Gamma_g/ \Gamma_g(l)}\Theta_{\cM,\tilde\cM}(\tau)^r|M
$$
where $\Theta_{\cM,\tilde\cM}(\tau)^r|M$ denotes
the usual slash operator and $M$ runs through a set of representatives of
$\Gamma_4/\Gamma_4(2l)$. Weissauer then shows that for given
$\varepsilon>0$ and $\tau\in\HH_4$ there is a subset $\tilde\cM
\subset \cM$ such that $\Theta_{\cM,\tilde\cM}(\tau)\neq 0$ and such
that the resulting form $F_r$ has the property that $m/k\ge
1/(12+\varepsilon)$.

We have to compare the vanishing order of such a form $F_r$ on
$\DVOR4$ with its vanishing order on~$E$. In order to do this, we
consider the $2$-dimensional cone $\Span{x^2_1, e}$. The dual cone was
computed in subsection~\ref{depth4section} above,
equation~\eqref{dualx1}. We write $m'=(m_1, m_2, m_3, m_4)$ and assume
that we have normalised in such a way that $-1/2\le m_i\le 1/2$. In
order to compute the vanishing order of $F_r(\tau)$ we have to compute
the Fourier expansions of the theta constants
$$
\Theta_{m, 0}(\tau, 0)=\sum_{q\in\ZZ^4} \prod_{i,j} 
t_{ij}^{\half(m_i+q_i)(m_j+q_j)} e^{2\pi i(q+m')^t {m''}}.
$$
We can rewrite this in terms of the coordinates $T_1, T_2,\ldots, T_{10}$. The
vanishing order of $\Theta_{m,0}(\tau, 0)$ along the divisor corresponding to
$x^2_1$ is then the minimum of the exponent of
$T_1^{\half (q_1+m_1)^2}$ for $q_1\in \ZZ$ and equals $m^2_1/2$.

The divisor $E$ is given by $T_2=0$, so the vanishing order of
$F_r(\tau)$ along~$E$ is the minimum over all $q\in\ZZ^4$ of the
exponents of $T_2$ in the summand given by~$q$. For fixed $q$ this
order can easily be computed to be $\half e(q+m)$, where $e$ is
the familiar quadratic form given in Section~\ref{toroidal},
equation~\eqref{definee}. That is,
\begin{eqnarray*}
\hf e(q+m)&=&\sum_i(q_i+m_i)^2+(q_1+m_1)(q_2+m_2)\\
&&\mbox{}-(q_1+m_1)(q_3+m_3)-(q_1+m_1)(q_4+m_4)\\
&&\mbox{}-(q_2+m_2)(q_3+m_3)-(q_2+m_2)(q_4+m_4).
\end{eqnarray*}
For $x=(x_1,x_2,x_3,x_4)\in\RR^4$ we define 
$$
\emin(x)=\min_{q\in\ZZ^4} e(q+x);
$$
then $F_r(\tau)$ will have the required vanishing as long as
$$
\sum_{m\in\cM\setminus\tilde\cM} \big(\emin(m)-\sfrac{9}{2}m^2_1\big) \ge 0.
$$
We claim that this is true for $l$ large enough. Given $\varepsilon>0$
we have $\#\tilde\cM<\varepsilon\#\cM$ for $l\gg 0$, so
(cf.~\cite[pp. 218--219]{Wei})
\begin{eqnarray*}
\lefteqn{\lim_{l\to\infty}\frac{1}{\#(\cM\setminus\tilde\cM)}
\sum_{m\in\cM\setminus\tilde\cM} \big(\emin(m)-\sfrac{9}{2}m^2_1\big)}\\
\ &=&
\int\!\int\!\int\!\int\nolimits_{[-\half,\half]^4}
\big(\emin(x)-\sfrac{9}{2}x^2_1\big)\,dx_1\,dx_2\,dx_3\,dx_4\\
&=&-\frac{3}{16}+\int\!\int\!\int\!\int\nolimits_{[-\half,\half]^4}
\emin(x)\,dx_1\,dx_2\,dx_3\,dx_4.
\end{eqnarray*}
The integral is not easy to evaluate, even though its value is
rational. The region of $\RR^4$ for which the minimum is achieved by
some particular value of~$q$ is a Delaunay cell for the quadratic
form~$e$, but these are complicated: there is a complete description
in~\cite{V2a}. Instead of attempting to evaluate the integral
precisely, we chose to estimate it by calculating
$\frac{1}{\#\cM'}\Sum_{m\in\cM'} \emin(m)$ using a computer, for a
suitable set of points~$\cM'$. Taking $\cM'$ to be the set of points
with coordinates of the form $\frac{2k+1}{158}$ gave us the estimate
\begin{equation*}
\frac{1}{\#\cM'}\Sum_{m\in\cM'} \emin(m)\approx 0.2166667
\end{equation*}
and by bounding the derivatives of the piecewise differentiable
continuous function $\emin$ one easily checks that the difference
between this and the actual value of the integral is less
than~$0.025$. Therefore $\int \emin > \frac{3}{16} = \int
\frac{9}{2}x^2_1$, and this proves the result.
\end{Proof}
\begin{remark}
The numerical evidence is overwhelming that $\int
\emin=\frac{13}{60}$, so that 
$$
\lim\limits_{l\to \infty}\frac{1}{\#(\cM\setminus\tilde\cM)}
\sum_{m\in\cM\setminus\tilde\cM}
\big(\emin(m)-\sfrac{9}{2}m^2_1\big) =\frac{7}{240}.
$$
\end{remark}

\subsection{Curves of depth 1}\label{depth1section}

In this subsection we want to prove a result (Proposition~\ref{nefondepth1})
giving
conditions for a divisor to have non-negative intersection with every
curve of depth~$1$. If $C\subset\AVOR4$ is a curve of depth~$1$ we
denote the
irreducible components of $\alpha_{n,\Vor}^{-1}(C)\subset\AVOR4(n)$
by $C_j(n)$.
Since $C$ is of depth~$1$ it is contained in $\DVOR4$: any component
$C_j(n)$ is therefore contained in a boundary component $\DVOR{4,i}(n)$.

The crucial point of the proof is the following lemma.
\begin{lemma}\label{goodJsection}
For any curve $C$ as above, given $\varepsilon >0$,
there exist integers $k$ and $n$ and a boundary component
$\DVOR{4,i}(n)$, for which we can find a section $\bds\in
H^0\big(kJ_4(n)\big)$ and a component $C_j(n)\subset \DVOR{4,i}$ such that
\begin{enumerate}
\item [(i)] $\bds|_{C_j(n)}\not\equiv 0$,
\item [(ii)] $\bds$ vanishes on $p_{i,n}^*\big(D_3(n)\big)$
to order $\lambda$, with $\lambda/k\ge n/(12+\varepsilon)$.
\end{enumerate}
\end{lemma}
\begin{Proof}
The proof of this lemma is a version of the argument of Weissauer
which we have already used for curves of depth~$0$. The argument given
in \cite[Proposition 4.1]{Hu} for $g=2$ is valid verbatim for all
$g\ge 2$.
\end{Proof}
\begin{proposition}\label{nefondepth1}
Let $C$ be a depth~$1$ curve and let $H=aL-b D_4-c E$ be a divisor
on $\AVOR4$ with $a\ge 0$, $a-12b\ge 0$ and $b\ge c\ge 0$. Then $H.C\ge 0$.
\end{proposition}
\begin{Proof}
It is sufficient to prove the result with the stronger condition
$a-12b>0$, since we can then take the limit as $a-12b\to
0$. Furthermore, because the cover $\alpha_{n,\Vor}$ is Galois, it is
enough to prove that $\alpha_{n,\Vor}^*(H).C_j(n)\ge 0$ for some $n$
and some component $C_j(n)$. We first of all choose some
$\varepsilon>0$ such that
$$
(a-12b)-b\left(1-\frac{12}{12+\varepsilon}\right)>0.
$$
We also choose a point $(z,\tau)\in\CC^3\times\HH_3$ whose image
$$
[(z,\tau)] \in D'_4=\big((\ZZ^3\times \ZZ^3) \rtimes \Sp(6,\ZZ)\big)
\backslash \CC^3\times \HH_3
$$
lies on the curve $C$. Choose a boundary component $D'_{4,1}(n)$ and consider
the point
$$
[(z,\tau)]_{n,1} \in D'_{4,1}(n)\imic\big((n\ZZ^3\times n\ZZ^3) \rtimes
\Gamma_3(n)\big)\backslash\CC^3\times \HH_3.
$$
It lies on some component $C_j(n)$ of the preimage of $C$ in
$\DVOR{4,1}(n)$. We shall prove that $\alpha_{n,\Vor}^*(H).C_j(n)>0$
for this component.

Recall from~\eqref{defineM} and~\eqref{defineJ} that
\begin{eqnarray*}
-\DVOR{4,1}(n)|_{\DVOR{4,1}(n)}=-\DVOR{4,1}(n)|_1&=&\sfrac 1 n M_{4,1}(n)
-\sfrac 1 n L\\
&=&\sfrac 1 n J_4(n)+E(n)|_1-\sfrac 1 n L,\\
\end{eqnarray*}
where $E(n)|_1=\sum_s E_s(n)|_1$. It follows (see Remark~\ref{Galois}
and Proposition~\ref{pullbackD3}) that
\begin{eqnarray*}
\alpha_{n,\Vor}^*(H)|_1&=&
aL-bn\DVOR{4,1}(n)|_1-bn\Sum_{i\neq 1} \DVOR{4,i}(n)|_1\\
&&\mbox{}-bnE(n)|_1-cnE(n)|_1\\
&=&(a-b)L+bJ_4(n)-bn p_{1,n}^*\big(bD_3(n)\big)+n(b-c)E(n)|_1.
\end{eqnarray*}

In terms of divisors, Lemma~\ref{goodJsection} means that
for some divisor $B\not\supset C$
$$
J_4(n)\sim\sfrac 1 k B+\sfrac{\lambda}{k}p_{i,n}^*\big(D_3(n)\big).
$$
Taking also into account that the divisor $L$ on $\DVOR{4,1}(n)$ is
a pullback from $\AVOR3(n)$ we find that
$$
\alpha_{n,\Vor}^*(H)|_1=p_{1,n}^* \left((a-b)L-b\Big(n-
\sfrac{\lambda}{k}\Big)D_3(n)\right) + \sfrac b k B+n(b-c)E(n)|_1.
$$
By construction $B.C_j(n)\ge 0$ and since $C$ is a curve of
depth~$1$ and $b\ge c$ we also have $n(b-c)E(n)|_1.C\ge 0$. 
The result now follows from the corresponding result, \cite[Theorem
0.2]{Hu}, on $\AVOR3(n)$, provided
$$
a-b>12\frac b n \left(n-\frac{\lambda}{k}\right).
$$
Since $\lambda k/n \ge 1/(12+\varepsilon)$ this follows from our choice of
$\varepsilon$.
\end{Proof}
\subsection{Curves of depth 2}\label{depth2section}

Let $C\subset \AVOR4$ be a curve of depth $\ge 2$ which is not
contained in the exceptional divisor $E$. Then (see
Remark~\ref{strata}) there are at least two different linear forms
$l\neq l'$ such that $C$ is contained in the divisors corresponding to
the rays $\Span{l^2}$ and $\Span{\smash{{l'}^2}\vphantom{l^2}}$. This
leads us to study the intersection of two boundary divisors in
$\AVOR4(n)$. Assume that $n\ge 3$ and that $\DVOR{4,ij}(n)\neq
\emptyset$. We have seen that each of the boundary components admits a
fibration
$$
p_i\colon\DVOR{4,i}(n)\To\AVOR3(n).
$$
Recall also (see e.g. \cite{Hu}) that each boundary component $D_{3,k}(n)$ of
$\AVOR3(n)$ admits a fibration
$$
q_k\colon D_{3,k}(n)\to\AVOR2(n).
$$
Indeed, this is the universal family over $\AVOR2(n)$. Since $\DVOR{4,ij}(n)$
is mapped to the boundary of $\AVOR3(n)$ this gives rise to the following
situation, for some $k=k(i,j)$ determined by the ordered pair~$(i,j)$:
$$
\CD
\DVOR{4,ij}(n)\  @.\subset \  \DVOR{4,i}(n)\\
@VVp_i|_j V @ VV p_i V\\
D_{3,k}(n)\  @. \subset \ \AVOR3(n)\\
@VV q_k V @. \\
\AVOR2(n).
\endCD
$$
Let $r_{ij}=q_k\circ p_i|_j$. By Proposition~\ref{pullbackD3}
we have that
\begin{eqnarray*}
-D_4(n)|_{ij}&=&-\DVOR{4,i}(n)|_{ij}-p_i|_j^* \big(D_3(n)\big)\\
&=&-\DVOR{4,i}{(n)}|_{ij}-p_i|_j^*\big(D_{3,k}(n)|_k\big)
-r_{ij}^*\big(D_2(n)\big)
\end{eqnarray*}
where $D_2(n)$ is the boundary of $\AVOR2(n)$, since
\begin{equation}\label{pullbackD2}
D_3(n)=D_{3,k}(n)+q^*_k\big(D_2(n)\big).
\end{equation}

\begin{lemma}\label{restrDij}
For $H=aL-bD_4(n)-cE(n)$ we have
\begin{eqnarray*}
H|_{ij}&=&\left(a-2\sfrac bn\right)L
+\sfrac b n p_i|_j^* M_{3,k}(n)
+\sfrac b n p_j|_i^* M_{3,k}(n)
-br_{ij}^*\big(D_2(n)\big)\\
&&\mbox{}
-br_{ji}^*\big(D_2(n)\big)
+b\Sum_{m\neq i,j}
\DVOR{4,m}|_{ij}+(4b-c)E(n)|_{ij}.
\end{eqnarray*}
\end{lemma}
\begin{Proof}
It follows from Proposition~\ref{pullbackD3} that
$$
p_i|_j^*\big(D_3(n)|_{k}\big)=
\Sum_{l\neq i}\DVOR{4,l}(n)|_{ij}+4E(n)|_{ij}
$$
and hence
$$
-\DVOR{4,j}(n)|_{ij}= -p_i|_j^*\big(D_3(n)|_{k}\big)
+\Sum_{l\neq i,j} \DVOR{4,l}(n)|_{ij}+ 4E(n)|_{ij}.
$$
Again using equation~\eqref{pullbackD2}, this implies
$$
-\DVOR{4,j}(n)|_{ij}= -p_i|_j^*\big(D_{3,k}(n)|_{k}\big)
-r^*_{ij}(D_2(n))+\Sum_{l\neq i,j} \DVOR{4,l}(n)|_{ij}+4E(n)|_{ij}.
$$
Using
$$
-D_{3,k}(n)|_{k}=\sfrac{1}{n} M_{3,k}(n)-\sfrac{1}{n}L
$$
we find that
\begin{eqnarray}\label{singlecpt}
-\DVOR{4,j}(n)|_{ij}&=& \sfrac{1}{n} p_i|_j^*\big(M_{3,k}(n)\big)
-\sfrac{1}{n}L-r_{ij}^*\big(D_2(n)\big)\nonumber\\
&&\mbox{}+\Sum_{l\neq i,j} \DVOR{4,l}(n)|_{ij}+4E(n)|_{ij}
\end{eqnarray}
Applying this formula also with $i$ and $j$ interchanged, and using
the fact that $-D_4(n)=-\Sum_l \DVOR{4,l}(n)-4E(n)$, we obtain
\begin{eqnarray*}
-D_4(n)|_{ij}&=&
\sfrac{1}{n}p_i|_j^*\big(M_{3,k}(n)\big)+
\sfrac{1}{n}p_j|_i^*\big(M_{3,k}(n)\big)-\sfrac{2}{n}L
-r_{ij}^*\big(D_2(n)\big)\\
&&\mbox{}
-r_{ji}^*\big(D_{2}(n)\big)
+\Sum_{m\neq i,j}\DVOR{4,m}(n)|_{ij} + 4E(n)|_{ij}.
\end{eqnarray*}
(There is a minor abuse of notation here, since the $k$s are not the same;
but this does not matter.) The result follows from this immediately.
\end{Proof}

We need to understand the last term in the expression in Lemma~\ref{restrDij}.
A component $E_s(n)|_{ij}$ of the restriction of $E(n)$ to
$\DVOR{4,ij}(n)$ corresponds to a $3$-dimensional Voronoi cell
$\Span{l_i^2, l_{\smash{j}}^2, e_s}\in\Vor(4)$ where $l_i$ and $l_j$ are
linear forms and $e_s$ is $\GL(\LL_4)$-equivalent to~$e$. These were
classified up to $\GL(\LL_4)$-equivalence in Corollary~\ref{dim3cones}.

Let $E$ and $E'$ be the components of the exceptional divisor
corresponding to the elements $e$ and $e'$ of $\MM_4$ of
Corollary~\ref{dim3cones}. As usual we assume that $\DVOR{4,1}(n)$ and
$\DVOR{4,2}(n)$ correspond to $x_1^2$ and $x_2^2$ respectively. Note
that the divisor $D_{3,k}(n)\subset\AVOR3(n)$, the image of $p_1|_2$,
corresponds to $x^2_2\in\MM_3$.
\begin{proposition}\label{fourcpts}
Let $E_s(n)|_i$ be a component of the restriction of $E(n)$ to some divisor
$\DVOR{4,i}(n)$. Then $p_i\big(E_s(n)|_i\big)$ is contained in exactly four
boundary components of $\AVOR3(n)$.
\end{proposition}

\begin{Proof}
Since we now have only one non-exceptional boundary component to
consider, it is enough to do this for the standard exceptional
divisor~$E$, corresponding to~$\eta$. This is because the stabiliser
of $\Pi_2(4)$ fixes $\eta$ and permutes the generators. But $\bar
e=\Sym_2\pr_1(e)$ is a sum of four forms of rank~$1$ in $\MM_3$, by
equation~\eqref{proje} in the proof of Proposition~\ref{exccones}.
These rank~$1$ forms together span a cone of $\Vor(3)$ and $p_i(E)$ is
contained in the intersection of the four corresponding boundary
components.
\end{Proof}

\begin{remark}\label{todeepestpoints}
In the two cases $E$ and $E'$ above we find two different kinds of
behaviour after projecting twice. The image $p_1(E|_1)$ is not
contained in $D_{3,k}(n)$, but $p_1(E'|_1)$ is. This is because
$$
\bar e = (x_2-x_3)^2+(x_2-x_4)^2+x_3^2+x_4^2
$$
does not have $x_2^2$ as a summand. In this case $r_{12}(E|_{12})$ is
contained in two boundary components of $\AVOR2(n)$.

On the other hand
$$
\bar e'= (x_2-x_3)^2+(x_3-x_4)^2+x_2^2+x_4^2
$$
does have $x_2^2$ as a summand, so $p_1(E'|_1)\subseteq D_{3,k}(n)$.
So $r_{12}(E'|_{12})$ is contained in three different boundary
components of $\AVOR2(n)$. In other words, $E'|_{12}$ is mapped under
$r_{12}$ to a deepest point in $\AVOR2(n)$, whereas $E|_{12}$ is mapped
to the intersection of two boundary components, i.e. to a $\PP^1$ in
$\AVOR2(n)$ if~$n\ge 3$.
\end{remark}

\begin{cor}\label{coeffs}
Let $E_s(n)|_{ij}$ be a component of $E(n)|_{ij}$. Then the coefficient of
$E_s(n)|_{ij}$ in $r_{ij}^*\big(D_2(n)\big)$ is equal to~$3$ if $E_s(n)$ is
mapped under $r_{ij}$ to a deepest point in $\AVOR2(n)$ and equal to~$4$
otherwise.
\end{cor}
\begin{Proof}
$D_2(n)$ is given by the support function $\psi_2$ on $\Vor(2)$ that takes the
value~$1$ on the primitive generator of every ray. In toric terms,
$q_k$ is given by the projection $\Sym_2\pr_2\colon\MM_3\to\MM_2$, which maps
$\bar e$ to $2(x_3^2+x_4^2)$ and $\bar e'$ to
$x_3^2+x_4^2+(x_3-x_4)^2$. So
$\psi_2(\bar e)=2\big(\psi_2(x_3^2)+\psi_2(x_4^2)\big)=4$ and similarly
$\psi_2(\bar e')=3$.
\end{Proof}

We are now in a position to begin checking nefness for depth~$2$ curves.
Let $C$ be such a curve. We choose a maximal cone $\sigma$
such that $C$ is contained in the closure of (the image in the moduli
space) of $\cO(\sigma)$. Since $C$ is a depth~$2$ curve, it is not
contained in the exceptional divisor and hence $\sigma$ must be of the form
$\Span{l^2_1,\ldots, l^2_k}$, where the $l_i$ are linear forms on $\LL_4$ and
where $\{l_1=\cdots=l_k=0\}$ is a plane in $\LL_4\otimes\RR$. In particular,
we may assume that $l_1$ and $l_2$ are linearly independent and that the other
linear forms $l_k$, $k\ge 3$ are linear combinations of $l_1$ and $l_2$.
Up to $\GL(\LL_4)$-equivalence we may assume that $\DVOR{4,i}(n)$ and
$\DVOR{4,j}(n)$ are $\DVOR{4,1}(n)$ and $\DVOR{4,2}(n)$, corresponding
to $\Span{l^2_1}=\Span{x^2_1}$ and $\Span{l^2_2}
=\Span{x^2_2}$ respectively. We may regard $\sigma$ as a cone in $\Vor(2)$
and from the known description of $\Vor(2)$ it follows that we need only
consider the cases $\sigma=\Span{x^2_1, x^2_2}$ or
$\sigma=\Span{x^2_1, x^2_2, (x_1-x_2)^2}$.

In particular $C$ is in either two or three (necessarily non-exceptional,
since $\depth(C)\neq 4$) irreducible components of the boundary
of~$\AVOR4(n)$. For any curve $C$ with $0<\depth(C)<4$ we define the
boundary multiplicity of $C$ to be the number $\mu(C)$ of irreducible
components of the boundary that contain~$C$: it
is the multiplicity of the generic point of $C$ as a point of $D_4(n)$.

First suppose that $\depth(C)=2$ and $\mu(C)=2$: that is, there are just
two such components,
which with the assumptions above are $\DVOR{4,1}$ and $\DVOR{4,2}$.
We shall use the expression for $H|_{ij}$ which we derived in Lemma
\ref{restrDij}. 

Notice (see Remark~\ref{todeepestpoints}) that if a component of
$E(n)|_{12}$ is contracted to a (deepest) point by $r_{12}$ then it is
also contracted to a point by $r_{21}$. Correspondingly we decompose 
$E(n)|_{12}$ as
$$
E(n)|_{12}=E_+(n)+E_-(n),
$$
where $E_+(n)$ consists of all the components that are not contracted
to points.

\begin{lemma}\label{M3versusD2}
Suppose $n\ge 3$ and $C\subset\DVOR{4,12}(n)$ is a depth~$2$ curve of
boundary multiplicity~$2$. Then, given $\varepsilon > 0$, we can write
$$
\frac 1 n p_1|_2^* \big(M_{3,k}(n)\big)\sim R_1
+\lambda_1 r_{12}^*\big(D_2(n)\big)+E_+(n)
$$
where $R_1$ is an effective $\QQ\,$-divisor with $C\not\subset \Supp R_1$
and $\lambda_1\ge 1/(12+\varepsilon)$. A similar statement holds for
$\frac 1 n p_2|_1^* \big(M_{3,k}(n)\big)$.
\end{lemma}
\begin{Proof}
We have already explained that one can construct suitable sections of
some power of $M_{3,k}(n)$ by taking products of theta functions of
the form $\Theta_{m' m''}(z,\tau)$, where
$$
\tau=\begin{pmatrix}\tau_{33} &
\tau_{34}\\\tau_{34} & \tau_{44}\end{pmatrix}\in
\HH_2\mbox{ and }z=(z_1, z_2)=(\tau_{23}, \tau_{24}).
$$
The claim about the vanishing
along $r_{12}^*\big(D_2(n)\big)$ follows by Weissauer's argument as in
\cite[Proposition 4.1]{Hu}. To check the contribution along the exceptional
divisors it is sufficient to check the divisors given by
the rays $\eta$: the terms corresponding to $\eta'=\Span{e'}$
(see Lemma~\ref{dim3cones}) can be absorbed into~$R_1$.
This works as in the proof of Proposition~\ref{sectionsofJ}.
The Fourier expansion, for $q=(q_3,q_4)\in\ZZ^2$, reads
\begin{eqnarray*}
\lefteqn{\Theta_{m' m''}(z,\tau)=}\\
&&\hskip-12pt \Sum_{q\in\ZZ^2}
t_{33}^{\half(q_3+m'_3)^2n} t_{44}^{\half
(q_4+m'_4)^2n}
t_{34}^{(q_3+m'_3)(q_4+m'_4)n} t_{23}^{(q_3+m'_3)n}
t_{24}^{(q_4+m'_4)n} e^{2\pi i(q+m')^t m''}.
\end{eqnarray*}
Let $w_3=(q_3+m'_3)$ and $w_4=(q_4+m'_4)$. A computation analogous to
that in the proof of Proposition~\ref{sectionsofJ} shows that the
vanishing order along $E_+(n)$ is equal to 
the minimum value of $2+2(w_3^2+w_4^2-w_3-w_4)$ for these values of
$w_3$ and~$w_4$.
The real function $w_3^2+w_4^2-w_3-w_4$ assumes its minimum
at $w_3=w_4=1/2$ where its value is $-1/2$ and this gives the result.
\end{Proof}
\begin{proposition}\label{nefondepth22}
Let $C$ be a depth~$2$ curve of boundary multiplicity~$2$,
and let $H=aL-b D_4(n)-c E(n)$ be a divisor
on $\AVOR4(n)$ with $a-12 b/n\ge 0$, $b\ge c\ge 0$. Then $H.C\ge 0.$
\end{proposition}
\begin{Proof}
Using Proposition~\ref{restrDij} and Lemma~\ref{M3versusD2} we find that
\begin{eqnarray*}
H|_{12}&=&(a-2\sfrac bn)L+b(R_1+R_2)
+b\lambda_1 r_{12}^*\big(D_2(n)\big)+b\lambda_2 r_{21}^*\big(D_2(n)\big)\\
&&\mbox{}+2bE_+(n)
-b r_{12}^*\big(D_2(n)\big)-b r_{21}^*\big(D_2(n)\big)\\
&&\mbox{}+b\Sum_{i\neq 1,2} \DVOR{4,i}|_{12}+(4b-c)E(n)|_{12}.
\end{eqnarray*}
We can rewrite this in the form
\begin{eqnarray*}
H|_{12}&=&\Sum_{i=1,2}r_{ii'}^*\Big(\left(\sfrac a 2-\sfrac b n\right)L
-b\left(\hf -\lambda_i\right) D_2(n)\Big)+b(R_1+R_2)\\
&&\mbox{}-\Sum_{i=1,2} \sfrac b 2 r_{ii'}^* \big(D_2(n)\big)+
b\Sum_{i\neq 1,2} \DVOR{4,i}(n)|_{12}\\
&&\mbox{}+2bE_+(n)+(4b-c)E(n)|_{12},
\end{eqnarray*}
where $i'=3-i$. As before (Proposition~\ref{nefondepth1}), we may
assume that in fact $a-12b/n>0$. The first two summands then have
non-negative intersection with~$C$. This follows by induction from our
knowledge of the nef cone of $\AVOR2(n)$ (\cite[Theorem 0.2]{Hu}) and
the inequality $\lambda_i\ge 1/(12+\varepsilon)$, where we can assume
$\varepsilon>0$ arbitrarily small.

By construction also $(R_1+R_2).C\ge 0$. By Corollary~\ref{coeffs} we have
$$
r_{12}^*\big(D_2(n)\big)=\Sum_i \DVOR{4,i}(n)|_{12}
+4E_+(n) +3E_-(n)
$$
where $\DVOR{4,i}(n)$ runs through all components such that
$r_{12}|_i$ is not dominant, and a similar formula
for~$r_{21}^*$. Altogether we see that the coefficients of
$\DVOR{4,i}(n)|_{12}$ for $i\neq 1$,~$2$ and those of $E_+(n)$ and $E_-(n)$
are all non-negative. 
Since $C$ is
not contained in any of these divisors we have proved the assertion.
\end{Proof}

We now move on to the case where the curve $C$ has depth~$2$ and
boundary multiplicity~$3$. Such a curve is contained in the closure of
the image in $\AVOR4(n)$ of $\cO(\sigma)$ where $\sigma$ has
dimension~$3$, and thus in $\DVOR{4,I}(n)$ for some set $I$ of three
indices. For convenience we take $I=\{1,2,3\}$, and we denote by
$\gothS_3$ the symmetric group on three elements acting as the
symmetry group of~$I$.

For each $\xi\in \gothS_3$ we define the bundles
$$
\cM_\xi=\left(p_{\xi(1)}|_{\xi(2)}^*\big(M_{3,k(\xi(1),\xi(2))}
\big)\right)|_{\xi(3)}
$$
and
$$
\cD_\xi=\left(r_{\xi(1)\xi(2)}^*\big(D_2(n)\big)\right)|_{\xi(3)}
$$
on~$\DVOR{4,I}(n)$.
\begin{lemma}\label{S3symmetrisation}
With the above notation
$$
-4D_4(n)|_I=\sfrac{1}{n}\Sum_\xi \cM_\xi
-\sfrac{6}{n} L-\Sum_\xi \cD_\xi
+2\Sum_{i\not\in I} \DVOR{4,i}(n)|_I+8E(n)|_I.
$$
\end{lemma}
\begin{Proof}
Apply equation~\eqref{singlecpt} in the proof of Lemma~\ref{restrDij} with
$i=\xi(1)$ and $j=\xi(2)$ and restrict to $\DVOR{4,\xi(3)}(n)$. Rearranging
this gives
$$
-\DVOR{4,\xi(2)}(n)|_I-\DVOR{4,\xi(3)}(n)|_I=
\sfrac{1}{n}\cM_\xi-\sfrac{1}{n} L - \cD_\xi
+\Sum_{i\not\in I} \DVOR{4,i}(n)|_I+4E(n)|_I;
$$
taking the sum over $\xi\in\gothS_3$ gives
$$
-4\Sum_{i\in I} \DVOR{4,i}(n)|_I=\sfrac{1}{n}\Sum_\xi \cM_\xi
-\sfrac{6}{n} L-\Sum_\xi \cD_\xi
+6\Sum_{i\not\in I} \DVOR{4,i}(n)|_I+24E(n)|_I.
$$
Since $D_4(n)|_I=\Sum_i\DVOR{4,i}(n)|_I+4E(n)|_I$ we have the formula stated.
\end{Proof}
In this case we again decompose $E(n)|_I$ as $E_+(n)+E_-(n)$ by
assigning a component to $E_-$ if it is contracted to a deepest point
by the~$r_{ij}$. Again a component is either contracted by all of the
$r_{ij}$ or by none of them. This is easy to see, since $\sigma$ is
$\GL(\LL_4)$-equivalent to $\Span{x_1^2,x_2^2,(x_1-x_2)^2}$ and hence
$r_{ij}$ corresponds to the linear map $\LL_4\to\LL_2$ with kernel
spanned by $x_1$ and $x_2$, independently of~$i$ and~$j$.
\begin{lemma}\label{cMversusD2}
For any $\xi\in\gothS_3$, given $\varepsilon>0$ we can write
$$
\sfrac{1}{n}\cM_\xi\sim R_\xi+\lambda_\xi\cD_\xi +E_+(n)
$$
where $R_\xi$ is an effective $\QQ\,$-divisor on $\DVOR{4,I}(n)$ such that
$C\not\subseteq \Supp R_\xi$, and $\lambda_\xi>1/(12+\varepsilon)$.
\end{lemma}
\begin{Proof}
Immediately from Lemma~\ref{M3versusD2}.
\end{Proof}
Note that, by Corollary~\ref{coeffs}, the coefficient of $E_+(n)$ in $\cD_\xi$
is equal to $4$ and the coefficient of $E_-(n)$ is equal to~$3$.
\begin{proposition}\label{nefondepth23}
Let $C$ be a curve of depth~$2$ and boundary multiplicity~$3$, and let
$H-aL-bD_4(n)-cE(n)$ be a divisor on $\AVOR4(n)$ such that $a-12b/n\ge
0$ and $b\ge 2c\ge 0$. Then $H.C\ge 0$.
\end{proposition}

\begin{Proof}
As usual we may assume $a-12b/n>0$. By Lemma~\ref{S3symmetrisation} and
Lemma~\ref{cMversusD2} we have
\begin{eqnarray*}
H|_I&=&\left(a-\sfrac{3b}{2n}\right)L+\sfrac{b}{4}\Sum_\xi R_\xi
-\sfrac{b}{4}\Sum_\xi (1-\lambda_\xi)\cD_\xi
+\sfrac{b}{2}\Sum_{i\not\in I}\DVOR{4,i}(n)|_I\\
&&\mbox{}+2bE(n)|_I+\sfrac{3b}{2}E_+(n)-cE(n)|_I\\
&=&\Sum_{\xi}r_{\xi(1)\xi(2)}^*\left(\left(\sfrac{a}{6}-\sfrac{b}{4n}
\right)L-\sfrac{b}{2}\left(\sfrac{1}{3}-\sfrac{\lambda_\xi}{2}
\right)D_2(n)\right)|_{\xi(3)}
-\sfrac{b}{12}\Sum_\xi \cD_\xi\\
&&\mbox{}+\sfrac{b}{4}\Sum_\xi R_\xi+\sfrac{b}{2}
\Sum_{i\not\in I}\DVOR{4,i}(n)|_I
+(2b-c)E_-(n)+(\sfrac{7}{2}b-c)E_+(n).
\end{eqnarray*}
Furthermore
$$
\cD_\xi=\Sum_{l\in J(\xi)}\DVOR{4,l}(n)|_I+4E_+(n)+3 E_-(n),
$$
where $l\in J(\xi)$ runs through all boundary components such that the image
$r_{\xi(1)\xi(2)} \big(\DVOR{4,l}(n)|_I\big)$ is contained in the boundary of
$\AVOR2(n)$. Note that $I\cap J(\xi)=\emptyset$.

Using this we get
\begin{eqnarray*}
H|_I&=&\Sum_{\xi}r_{\xi(1)\xi(2)}^*
\left(\left(\sfrac{a}{6}-\sfrac{b}{4n}\right)L
-\sfrac{b}{2}\left(\sfrac{1}{3}-\sfrac{\lambda_\xi}{2}\right)
D_2(n)\right)|_{\xi(3)}+\sfrac{b}{4}\Sum_\xi R_\xi\\
&&\mbox{}+\Sum_{i\not\in I} d_i\DVOR{4,i}(n)|_I
+(\sfrac{b}{2}-c)E_-(n)+(\sfrac{3b}{2}-c)E_+(n).
\end{eqnarray*}
Moreover the coefficients $d_i$ are non-negative. By induction, i.e.\ 
by our knowledge of the nef cone of~$\AVOR2(n)$ from~\cite{Hu}, we can
assume that the divisor $\Big(\sfrac{a}{6}-\sfrac{b}{4n}\Big)L
-b\Big(\sfrac{1}{6}-\sfrac{\lambda_\xi}{4}\Big)D_2(n)$ is nef and
therefore that its pullback has non-negative intersection with~$C$.
Also $\sum_{\xi\in\gothS_3}R_\xi.C\ge 0$. Since $C$ is not contained
in any of the divisors $\DVOR{4,i}(n)$ for $i\not\in I$, nor in any
component of $E(n)|_I$, and since $d_i\ge 0$ and the coefficients of
the components of $E_-(n)$ and $E_+(n)$ are all non-negative, the
result follows.
\end{Proof}

\subsection{Curves of depth 3}\label{depth3section}

Supppose $C$ is an irreducible depth~$3$ curve of boundary
multiplicity $\mu=\mu(C)$. Then $3\le \mu\le 6$ and $C$ is contained
in $\DVOR{4,i}(n)$ if and only if $i\in I$ for some index set $I$ of
size~$\mu(C)$. We assume without loss of generality that
$I=\{1,2,\ldots,\mu\}$, and denote by $\gothS_\mu$ the symmetric group
on $\mu$ symbols acting as the symmetry group of~$I$. For any
$\xi\in\gothS_\mu$ we have the following diagram:
$$
\CD
\DVOR{4,I}(n)\  @.\subset \  \DVOR{4,\xi(1)\xi(2)}(n)
\ @.\subset \  \DVOR{4,\xi(1)}(n)\\
@VV p_{\xi(1)}|_I V @ VV p_{\xi(1)}|_{\xi(2)} V @ VV p_{\xi(1)} V\\
D_{3,K(\xi)}(n)\  @. \subset D_{3,k(\xi)}(n)\ @. \subset \ \AVOR3(n)\\
@VV q_{k(\xi)}|_{K(\xi)} V @ VV q_{k(\xi)} V\\
D_{2,{m(k(\xi))}}(n)\ @. \subset \ \AVOR2(n)\\
@VV s_{m(k(\xi))} V @. \\
\AVOR1(n).
\endCD
$$
Here the index $k(\xi)$ and the set of indices $K(\xi)$ are all
determined by the choice of $\xi\in\gothS_\mu$: we define $K(\xi)$ to
be the set of indices~$k$ such that
$p_{\xi(1)}|_I\big(\DVOR{4,I}(n)\big)\subset D_{3,k}(n)$. In fact
$k(\xi)$ depends only on $\xi(1)$ and $\xi(2)$, and $K(\xi)$ only on
$\xi(1)$. For any $k\in K(\xi)$, we define $m(k)$ so that
$q_k|_{K(\xi)}$ maps $D_{3,K(\xi)}(n)$ to $D_{2,m(k)}(n)$. It is not
quite immediate that $m(k)$ is well defined: in principle the image of
$q_k|_{K(\xi)}$ could be in the intersection of two boundary
components of $\AVOR2(n)$. However, this does not happen: see
Corollary~\ref{uniquem} below.

To reduce the amount of notation we write
$p_\xi=p_{\xi(1)}|_{\xi(2)}$, $r_\xi=r_{\xi(1)\xi(2)}
=q_{k(\xi)}\circ p_{\xi(1)}|_{\xi(2)}$, and
$r_\xi|_I=q_{k(\xi)}|_{K(\xi)}\circ p_{\xi(1)}|_I$. Note that $p_\xi$
and $r_\xi$ depend only on $\xi(1)$ and $\xi(2)$. As in the case of
depth~$2$ and boundary multiplicity~$3$ we define bundles
$$
\cD_\xi=\left(r_\xi^*\big(D_2(n)\big)\right)|_I
$$
on~$\DVOR{4,I}(n)$. We also write
$$
\cS_\xi=\left(\pxi\big(D_{3,k(\xi)}(n)\big)\right)|_I
$$
\begin{lemma}\label{Smusymmetrisation}
If $H=aL-bD_4(n)-cE(n)$ and $a\ge 12b/n\ge 0$, then there is a nef
$\QQ\,$-divisor $A_1$ such that
\begin{eqnarray*}
H|_I&=&A_1-\sfrac{b}{(\mu-1)\mu!}\Sum_{\xi\in \gothS_\mu}\cD_\xi\\
&&\mbox{}-\sfrac{b}{(\mu-1)(\mu-1)!}\Sum_{\xi\in \gothS_\mu}\cS_\xi
+\sfrac{b}{(\mu-1)}\Sum_{i\not\in I}\DVOR{4,i}(n)|_I
+\left(\sfrac{4b}{\mu-1}-c\right)E(n)|_I.
\end{eqnarray*}
\end{lemma}
\begin{Proof}
Replacing $\gothS_3$ by $\gothS_\mu$ in the proof of
Lemma~\ref{S3symmetrisation}, we see that
\begin{eqnarray*}
\lefteqn{-(\mu-1)(\mu-1)!\,D_4(n)|_I=}\\
&\;&-\Sum_{\xi\in\gothS_\mu} \cS_\xi
-\Sum_{\xi\in\gothS_\mu} \cD_\xi
+(\mu-1)!\Sum_{i\not\in I} \DVOR{4,i}(n)|_I
+4(\mu-1)!\,E(n)|_I.
\end{eqnarray*}
{}From this it follows that
\begin{eqnarray*}
H|_I&=&aL
-\sfrac{b}{(\mu-1)(\mu-1)!}\Sum_{\xi\in \gothS_\mu}\cS_\xi
-\sfrac{b}{(\mu-1)(\mu-1)!}\Sum_{\xi\in \gothS_\mu}\cD_\xi\\
&&\mbox{}+\sfrac{b}{(\mu-1)}\Sum_{i\not\in I}\DVOR{4,i}(n)|_I
+\left(\sfrac{4b}{\mu-1}-c\right)E(n)|_I\\
&=&\Sum_{\xi\in\gothS_\mu}r_\xi^*\left(
\sfrac{a}{\mu!}L-\sfrac{b}{\mu!}D_2(n)\right)|_I
-\sfrac{b}{(\mu-1)\mu!}\Sum_{\xi\in \gothS_\mu}\cD_\xi\\
&&\mbox{}-\sfrac{b}{(\mu-1)(\mu-1)!}\Sum_{\xi\in \gothS_\mu}\cS_\xi
+\sfrac{b}{(\mu-1)}\Sum_{i\not\in I}\DVOR{4,i}(n)|_I
+\left(\sfrac{4b}{\mu-1}-c\right)E(n)|_I,
\end{eqnarray*}
which gives the formula required. The first term is nef because it is
the pullback of a nef bundle from~$\AVOR2$.
\end{Proof}

For any boundary component $D_{3,k}(n)$ of $\AVOR3(n)$ we have as
before (equation~\eqref{pullbackD2})
\begin{equation}\label{pullbackD2(3)}
q^*_k \big(D_2(n)\big)|_k=\Sum_{j\neq k}D_{3,j}(n)|_k,
\end{equation}
so
$$
r_\xi^*\big(D_2(n)\big)|_I
=\Sum_{j\neq k(\xi)}p_{\xi(1)}^*\big(D_{3,j}(n)\big)|_I.
$$

\begin{lemma}\label{sizeofKxi}
For any $\xi\in\gothS_\mu$ (with $3\le \mu\le 6$) we have $\#K(\xi)=2$
or $\#K(\xi)=3$. If $\mu=3$ then $\#K(\xi)=2$; if $\mu=6$ then
$\#K(\xi)=3$. If $\mu=4$ or $\mu=5$ then both cases can occur,
depending on~$\sigma$ and~$\xi$.
\end{lemma}
\begin{Proof}
The cone $\sigma$ is of the form $\Span{l_1^2,\ldots,l_\mu^2}$ and
has three special properties: it is a Voronoi cone, $\sigma\in\Vor(4)$;
it contains no rank~$4$ forms (otherwise $\depth(C)=1$); and the
linear span of the~$l_i$ is of dimension~$3$. Then
$$
K(\xi)=\{\Sym_2\pr_{\xi(1)}(l_i)\mid \sigma=\Span{l_1,\ldots,l_\mu},\
i\neq\xi(1)\}.
$$
By using the $\GL(\LL_4)$-action we may assume that
$\sigma\prec\Pi_1(4)$, that is, that the generators of $\sigma$ are
$x_i^2$ or $(x_i-x_j)^2$; and we may assume that $x_4$ does not occur
at all, so that $\sigma\prec\Pi_1(3)$.

If $\mu=3$ then the $l_i$ are linearly independent, so the projections
of $l_{\xi(2)}$ and $l_{\xi(3)}$ are also linearly independent. So $\#
K(\xi)=2$.

If $\mu=4$, there are two possibilities. If three of the $l_i$ are
linearly dependent then $\# K(\xi)=2$ if $\xi(1)$ is one of those
three, since the projection identifies the other two; but $\#
K(\xi)=3$ if $\xi(1)$ is the fourth generator. On the other hand, if
no three of the $l_i$ are linearly dependent, then any two projections
are distinct, so~$\# K(\xi)=3$.

We remark that both these cases occur: examples are
$\Span{x_1^2,x_2^2, (x_1-x_2)^2, x_3^2}$ and $\Span{x_1^2,(x_1-x_2)^2,
x_3^2, (x_2-x_3)^2}$.

If $\mu=5$ then $\sigma$ is a codimension~$1$ face of $\Pi_1(3)$ and
thus equivalent to $\Span{x_1^2,x_2^2,x_3^2, (x_1-x_2)^2,
(x_1-x_3)^2}$. There are two linear relations involving three of the
$l_i$ and one generator (namely $l_1^2=x_1^2$) occurs in both of
them. Each linear relation involving $l_{\xi(1)}$ reduces $\#K(\xi)$
by~$1$, starting from $\mu-1$ (i.e.\ if there were no relations all
the other generators would give different elements of $K(\xi)$, so we
should have $\#K(\xi)=\mu-1$); so $\# K(\xi)=2$ if $\xi(1)=1$ and
$\#K(\xi)=3$ otherwise.

If $\mu=6$ then $\sigma=\Pi_1(3)$ and because $\GL(\LL_3)$ permutes
the generators we have $\#K(\xi)=\#K(\id)=3$ for all~$\xi$.
\end{Proof}

\begin{cor}\label{uniquem}
If $k\in K(\xi)$ then $m(k)$ is unique: that is, the image of
$q_k|_{K(\xi)}$ is contained in exactly one boundary component
$D_{2,m(k)}$ of~$\AVOR2(n)$.
\end{cor}
\begin{Proof}
The linear span of $K(\xi)$ has dimension~$2$, since it is the
projection of the linear span of the~$l_i$. Projecting again from an
element of~$K_(\xi)$ therefore gives a space of dimension~$1$.
\end{Proof}

We can characterise $m(k)$ by saying that it gives the
unique boundary component that contains the image of~$C$.
We denote the elements of $K(\xi)$ by $k(\xi)$, $k'(\xi)$ and (if $\#K(\xi)=3$)
$k''(\xi)$.

We also know from the case $g=2$ that $-D_{2,m(k)}(n)|_{m(k)}$ is
nef. We write
$$
D_2(n)=D_{2,m(k)}(n)+D_{2,\hat m(k)}(n),
$$
so $D_{2,\hat m(k)}(n)$ is the union of all the boundary components except
for the unique one that contains the image of~$C$.

\begin{lemma}\label{differentD3}
Let $D_{3,1}(n)$ and $D_{3,2}(n)$ be two boundary components of
$\DVOR3(n)$. Suppose $q_1\colon D_{3,1}(n)\to\AVOR2(n)$ and
$q_2\colon D_{3,2}(n)\to\AVOR2(n)$ are the associated projection maps,
and write $D_{2,m(1)}$ for the
image $q_1\big(D_{3,12}(n)\big)$ and similarly
$D_{2,m(2)}=q_2\big(D_{3,12}(n)\big)$. Then
\begin{enumerate}
\item [(i)] $q_2^* D_{2,\hat m(2)}(n)=q_1^* D_{2,\hat m(1)}(n)$,
\item [(ii)] $D_{3,1}(n)|_{12}=D_{3,2}(n)|_{12}
+q_2^* D_{2,m(2)}(n)|_{12}-q_1^* D_{2,m(1)}(n)|_{12}$.
\end{enumerate}
\end{lemma}
\begin{Proof}
\noindent(i) We may assume that $D_{3,1}(n)$ and $D_{3,2}(n)$ correspond to
$\Span{x_1^2}$ and $\Span{x_2^2}$ in $\Vor(3)$. Then $q_1$ is given by
$\Sym_2\pr_1$ so boundary components of the image of $q_1$ (which is
abstractly $\AVOR2(n)$) may be thought of as cones in the Voronoi
decomposition of quadratic forms in the variables $x_2$ and $x_3$,
while for $q_2$ one should consider quadratic forms in $x_1$ and
$x_3$. In particular $\hat m(1)$ is $\Span{x_2^2}$ and $\hat m(2)$ is
$\Span{x_1^2}$. Conside an arbitrary boundary component given by
$\Span{(a_1x_1+a_2x_2+a_3x_3)^2}\in\Vor(3)$. Under $q_1$ it maps to
$\Span{(a_2x_2+a_3x_3)^2}$, which is different from
$m(1)=\Span{x_2^2}$ if and only if $a_3\neq 0$. Similarly the image
under $q_2$ is different from $m(2)$ if and only if $a_3\neq 0$.

\noindent(ii) Applying (i) and using the equation (compare
equation~\eqref{pullbackD2})
\begin{eqnarray*}
-D_{3,2}(n)&=&\Sum_{i\neq 1,2}D_{3,i}(n)|_{12}-q_1^*D_2(n)|_{12}\\
&=&\Sum_{i\neq 1,2}D_{3,i}(n)|_{12}-q_1^*D_{2,m(1)}(n)|_{12}
-q_1^*D_{\hat m(1)}(n)|_{12}
\end{eqnarray*}
and the same equation with the indices interchanged, we obtain
\begin{eqnarray*}
\big(-D_{3,2}(n)+q_1^*D_{2,m(1)}\big)|_{12}
&=&\Sum_{i\neq 1,2}D_{3,i}(n)|_{12}-q_1^*D_{2,\hat m(1)}|_{12}\\
&=&\Sum_{i\neq 1,2}D_{3,i}(n)|_{12}-q_2^*D_{2,\hat m(2)}|_{12}\\
&=&\big(-D_{3,1}(n)+q_2^*D_{2,m(2)}\big)|_{12}
\end{eqnarray*}
as required.
\end{Proof}
\begin{proposition}\label{HwithoutSxi}
There is a nef $\QQ\,$-divisor $A_2$ such that the following expression
for $H|_I$ holds:
\begin{eqnarray*}
H|_I&=&A_2
-\sfrac{b}{(\mu-1)\mu!}\sum_{\xi\in\gothS_\mu}
r_\xi^* D_{2,\hat m(k)}(n)|_I\\
&&\mbox{}-\sum_{\xi\in\gothS_\mu}\sfrac{b}{(\#K(\xi)-1)(\mu-1)(\mu-1)!}
r_\xi^* D_{2,\hat m(k)}(n)|_I\\
&&\mbox{}+\sum_{\xi\in\gothS_\mu}\sfrac{b}{(\#K(\xi)-1)(\mu-1)(\mu-1)!}
\pxi \bigg(\sum_{l\not\in K(\xi)}D_{3,l}(n)\bigg)\\
&&\mbox{}+\sfrac{b}{\mu-1}\sum_{i\not\in I} \DVOR{4,i}(n)|_I
+\left(\sfrac{4b}{\mu-1}-c\right) E(n)|_I.
\end{eqnarray*}
\end{proposition}
\begin{Proof}
{}From equation~\eqref{pullbackD2(3)} we have
\begin{equation}\label{K2pbD2}
-D_{3,k'}|_{K(\xi)}=
\sum_{l\not\in K(\xi)}D_{3,l}|_{K(\xi)}-q_k^* D_2(n)|_{K(\xi)}
\end{equation}
if $\# K=2$ and 
\begin{equation}\label{K3pbD2}
-D_{3,k'}|_{K(\xi)}-D_{3,k''}|_{K(\xi)}=
\sum_{l\not\in K(\xi)}D_{3,l}|_{K(\xi)}-q_k^* D_2(n)|_{K(\xi)}
\end{equation}
if $\#K=3$. We also have, by Lemma~\ref{differentD3}
\begin{equation}\label{k'eqn}
D_{3,k}|_{K(\xi)}=D_{3,k'}|_{K(\xi)}+q_{k'}^*
D_{2,m(k')}(n)|_{K(\xi)}-q_k^* D_{2,m(k)}(n)|_{K(\xi)}
\end{equation}
and similarly for $k''$ in place of $k'$ if $\#K=3$. In the latter
case we add the two equations to obtain
\begin{eqnarray}\label{k''eqn}
D_{3,k}|_{K(\xi)}&=&\hf D_{3,k'}|_{K(\xi)}+\hf D_{3,k''}|_{K(\xi)} 
+\hf q_{k'}^* D_{2,m(k')}(n)|_{K(\xi)}\nonumber\\
&&\mbox{}+\hf q_{k''}^* D_{2,m(k'')}(n)|_{K(\xi)}
-q_k^* D_{2,m(k)}(n)|_{K(\xi)}.
\end{eqnarray}
We use these equations to eliminate $\cS_\xi$ from the formula in
Lemma~\ref{Smusymmetrisation}. {}From \eqref{K2pbD2} and \eqref{k'eqn}
we have
\begin{equation*}
\cS_\xi=-\sum_{l\not\in K(\xi)}\pxi D_{3,l}(n)|_{K(\xi)}+\cD_\xi
+\pxi q_{k'}^* D_{2,m(k')} 
- \pxi q_k^* D_{2,m(k)}
\end{equation*}
if $\#K=2$, and from \eqref{K3pbD2} and \eqref{k''eqn} we have
\begin{eqnarray*}
\cS_\xi&=&-\hf \sum_{l\not\in K(\xi)}\pxi D_{3,l}(n)|_{K(\xi)}
+\hf \cD_\xi\\
&&\mbox{}+\pxi\left(\hf q_{k'}^* D_{2,m(k')}
+\hf q_{k''}^* D_{2,m(k'')}
- q_k^* D_{2,m(k)}\right).
\end{eqnarray*}
The term $\frac{b}{\mu-1}\sum_{i\not\in I} \DVOR{4,i}(n)|_I
+\left(\frac{4b}{\mu-1}-c\right) E(n)|_I$ plays no role at this
point and we temporarily denote it by~$\ast$. So for $\#K=2$ we have
\begin{eqnarray}\label{K2}
H|_I&=&A_1+\ast-\sfrac{b}{(\mu-1)\mu!}\Sum_{\xi\in \gothS_\mu}\cD_\xi
-\sfrac{b}{(\mu-1)(\mu-1)!}\Sum_{\xi\in\gothS_\mu}\cD_\xi\nonumber\\
&&\mbox{}+\sfrac{b}{(\mu-1)(\mu-1)!}\Sum_{\xi\in \gothS_\mu} 
\sum_{l\not\in K(\xi)}\pxi D_{3,l}(n)|_{K(\xi)}\\
&&\mbox{}+\Sum_{\xi\in \gothS_\mu}
\pxi\left( q_{k'}^* D_{2,m(k')}
- q_k^* D_{2,m(k)}\right)\nonumber
\end{eqnarray}
and for $\#K=3$
\begin{eqnarray}\label{K3}
H|_I&=&A_1+\ast-\sfrac{b}{(\mu-1)\mu!}\Sum_{\xi\in \gothS_\mu}\cD_\xi
-\sfrac{b}{(\mu-1)(\mu-1)!}\Sum_{\xi\in\gothS_\mu}\hf\cD_\xi\nonumber\\
&&\mbox{}+\sfrac{b}{(\mu-1)(\mu-1)!}\Sum_{\xi\in \gothS_\mu} 
\sum_{l\not\in K(\xi)}\hf \pxi D_{3,l}(n)|_{K(\xi)}\\
&&\mbox{}-\sfrac{b}{(\mu-1)(\mu-1)!}\Sum_{\xi\in \gothS_\mu}
\pxi \left(\hf q_{k'}^* D_{2,m(k')}
+\hf q_{k''}^* D_{2,m(k'')}
- q_k^* D_{2,m(k)}\right).\nonumber
\end{eqnarray}
Since $-D_{2,m(k)}(n)|_{m(k)}$ is nef we may add those terms to
$A_1$, obtaining a nef $\QQ\,$-divisor $A_2$ where
\begin{equation*}
A_2=A_1-\sfrac{b}{(\mu-1)(\mu-1)!}\Sum_{\xi\in\gothS_\mu}
r_\xi^*D_{2,m(k(\xi))}(n)|_{m(k(\xi))}. 
\end{equation*}
In view of Lemma~\ref{differentD3} this
allows us to replace $\cD_\xi$ by $r_\xi^* D_{2,\hat m(k)}(n)|_I$.
This gives the desired coefficients for the $D_{2,\hat m(k)}(n)$ and
$D_{3,l}(n)$ terms. Finally, the terms
$$
\Sum_{\xi\in \gothS_\mu}
\left(\pxi q_{k'}^* D_{2,m(k')}
-\pxi q_k^* D_{2,m(k)}\right)
$$
and
$$
\Sum_{\xi\in \gothS_\mu}
\left(\hf \pxi q_{k'}^* D_{2,m(k')}
+\hf \pxi q_{k''}^* D_{2,m(k'')}
-\pxi q_k^* D_{2,m(k)}\right)
$$
vanish. Indeed, if we fix $\xi(1)=a$ say, and take
$\gothS_{\mu-1}=\{\xi\in\gothS_\mu \mid \xi(1)=a\}$ then
$$
\Sum_{\xi\in \gothS_{\mu-1}}
\left(p_a|_{\xi(2)}^*q_{k'}^* D_{2,m(k')}
-p_a|_{\xi(2)}^*q_k^* D_{2,m(k)}\right)=0
$$
and similarly for~$\#K=3$. 
\end{Proof}

In fact only the part of $H|_I$ supported on the exceptional divisors
could have negative intersection with~$C$. We have
\begin{equation}\label{coeffsgamma}
\pxi \bigg(\sum_{l\not\in K(\xi)}D_{3,l}(n)\bigg)|_I=
\sum_{j\in L(\xi)}\DVOR{4,j}(n)|_I+\sum_s \gamma_s(\xi) E_s(n)
\end{equation}
and
\begin{equation}\label{coeffsdelta}
r_\xi^* D_{2,\hat m(k)}(n)|_I=\sum_{j\in M(\xi)}\DVOR{4,j}(n)|_I+\sum_s
\delta_s(\xi) E_s(n),
\end{equation}
where $L(\xi)\supseteq M(\xi)$, $L(\xi)\cap I=\emptyset$ and
$\gamma_s$ and $\delta_s$ depend on the cone $\sigma$, the
permutation~$\xi$ and the component $E_s(n)$ of the exceptional
divisor.

\begin{cor}\label{gammadeltacor}
$H.C\ge 0$ provided that
\begin{eqnarray*}
0&\le&\sum_{\xi\in\gothS_\mu}\left(-\sfrac{b}{(\mu-1)\mu!}\delta_s(\xi)
-\sfrac{b}{(\#K(\xi)-1)(\mu-1)(\mu-1)!}\delta_s(\xi)\right.\\
&&\left.+\sfrac{b}{(\#K(\xi)-1)(\mu-1)(\mu-1)!}\gamma_s(\xi)\right)
+\sfrac{4b}{\mu-1}-c
\end{eqnarray*}
for every component $E_s(n)$ such that $E_s(n)\cap C\neq\emptyset$.
\end{cor}
\begin{Proof}
Away from the exceptional divisors $E_s(n)$, by
Proposition~\ref{HwithoutSxi} and equations~\eqref{coeffsgamma} and
\eqref{coeffsdelta} we can write
$$
-\sfrac{b}{(\mu-1)\mu!}\sum_{\xi\in\gothS_\mu}
r_\xi^* D_{2,\hat m(k)}(n)|_I +\sfrac{b}{\mu-1}\sum_{i\not\in I}
\DVOR{4,i}(n)|_I = \sum_{i\not\in I}
\alpha_i\DVOR{4,i}(n)|_I
$$
and
\begin{eqnarray*}
\lefteqn{-\sum_{\xi\in\gothS_\mu}\sfrac{b}{(\#K(\xi)-1)(\mu-1)(\mu-1)!}
r_\xi^* D_{2,\hat m(k)}(n)|_I}\\
&&\mbox{}+\sum_{\xi\in\gothS_\mu}\sfrac{b}{(\#K(\xi)-1)(\mu-1)(\mu-1)!}
\pxi \bigg(\sum_{l\not\in K(\xi)}D_{3,l}(n)\bigg)
=\sum_{i\not\in I} \beta_i\DVOR{4,i}(n)|_I
\end{eqnarray*}
with $\alpha_i\ge 0$ and $\beta_i\ge 0$. Hence
\begin{eqnarray*}
H|_I&=&A_2+\sum_{i\not\in I} (\alpha_i+\beta_i)\DVOR{4,i}(n)|_I\\
&&\mbox{}+\sum_s \bigg(\sum_{\xi\in\gothS_\mu}\bigg(-\sfrac{b}{(\mu-1)\mu!}
\delta_s(\xi)
-\sfrac{b}{(\#K(\xi)-1)(\mu-1)(\mu-1)!}\delta_s(\xi)\\
&&+\sfrac{b}{(\#K(\xi)-1)(\mu-1)(\mu-1)!}\gamma_s(\xi)\bigg)
+\sfrac{4b}{\mu-1}-c\bigg)E_s(n).
\end{eqnarray*}
Since $A_2$ is nef and $\DVOR{4,i}.C\ge 0$ for $i\not\in I$ we have
the result claimed.
\end{Proof}

We shall fix $s$ so that $E_s=E$, the component corresponding to the
ray~$\eta$, and drop the suffix~$s$ from the notation~$\gamma_s(\xi)$,
$\delta_s(\xi)$. Since we may assume that $C\cap E\neq \emptyset$, we
can take $\sigma$ to be a face of~$\Pi_2(4)$. Obviously we need only
consider the faces up to $G$-equivalence.

Unfortunately the inequality of Corollary~\ref{gammadeltacor} does not
always hold. We shall need an extra argument, applied in
Proposition~\ref{mu3special} and Proposition~\ref{mu6}, to handle the
cases where it fails.

We begin by calculating the values of $\gamma(\xi)$ and $\delta(\xi)$ for
three of the four cases of Lemma~\ref{dim3faces}. We shall work with
the representatives given in Proposition~\ref{dim3faces}, and with the
generators in the order given there. Note that the values of $\gamma(\xi)$ and
$\delta(\xi)$ do depend on~$\xi$, not just $\sigma$.

\begin{proposition}\label{gammadeltavalues}
If $\sigma$ is a $3$-dimensional face of $\Pi_2(4)$, not of type
RT\/${}^*$, and $\xi\in\gothS_3$ then the values of $\gamma(\xi)$ and
$\delta(\xi)$ are:
\begin{eqnarray*}
\gamma_(\xi)=3,\ \delta(\xi)=2&\quad&\mbox{ if }\sigma\mbox{ is of type
string and $\xi(1)\neq 3$;}\\
\gamma_(\xi)=2,\ \delta(\xi)=2&\quad&\mbox{ if }\sigma\mbox{ is of type
string and $\xi(1)= 3$;}\\
\gamma_(\xi)=2,\ \delta(\xi)=2&\quad&\mbox{ if }\sigma\mbox{ is of type
BF\/}^*;\\
\gamma_(\xi)=4,\ \delta(\xi)=4&\quad&\mbox{ if }\sigma\mbox{ is of type
disconnected.}
\end{eqnarray*}
\end{proposition}
\begin{Proof}
We consider the support functions $\psi_\gamma$ and $\psi_\delta$ on
(suitable copies of) $\Vor(3)$ and $\Vor(2)$ respectively that
determine the divisors $\sum_{l\not\in K(\xi)}D_{3,l}(n)$ and
$D_{2,\hat m(k)}(n)$: namely, $\psi_\gamma$ (respectively
$\psi_\delta$) takes the value~$0$ on a generator~$v$ of a ray in
$\Vor(3)$ (respectively $\Vor(2)$) if $v\in K(\xi)$ (respectively $v=
m(\xi)$) and~$1$ otherwise. Observe in particular that $\gamma(\xi)$
depends only on $\xi(1)$, since that determines~$K(\xi)$, but
$\delta(\xi)$ depends on the unordered pair $\{\xi(1),\xi(2)\}$.

String type: $\Span{x_1^2,x_2^2, x_3^2}$. Since $x_1^2$ and $x_2^2$
are interchanged by $k_3\in G$, which preserves~$\sigma$, we need only
consider the cases $\xi(1)=1$ and $\xi(1)=3$. If $\xi(1)=1$ then the
projection gives, as in~\eqref{proje}
\begin{equation}\label{projerepeated}
{\bar e}=(x_2-x_3)^2+(x_2-x_4)^2+x_3^2+x_4^2.
\end{equation}
In this case $K(\xi)=\{x_2^2,x_3^2\}$ so $\gamma(\xi)=3$. Taking the
further projection from $\xi(2)$ we get either ${\bar{\bar e}}=2x_3^2+2x_4^2$,
if $\xi(2)=2$, or ${\bar{\bar e}}=x_2^2+(x_2-x_4)^2+x_4^2$, if
$\xi(2)=3$. In the first case $m(\xi)=x_3^2$ and in the second case
$m(\xi)=x_2^2$: in either case $\delta(\xi)=\psi_\delta(\bar{\bar
e})=2$.

On the other hand, if $\xi(1)=3$ then the projection $\Sym_2\pr_3$
gives
$$
{\bar e}=x_1^2+x_2^2+x_4^2+(x_1+x_2-x_4)^2.
$$
$K(\xi)=\{x_1^2,x_2^2\}$ so $\gamma(\xi)=\psi_\gamma(\bar
e)=2$. Applying $k_1$ if necessary we may take $\xi(2)=1$, so
$\bar{\bar e}=x_2^2+x_4^2+(x_2-x_4)^2$ and $m(\xi)=x_2^2$ so
$\delta(\xi)=2$

BF${}^*$ type: $\Span{x_1^2,x_3^2, x_4^2}$. In this case the subgroup
of $G$ that preserves $\sigma$ acts transitively on the generators
($s_{23}k_2k_3$ interchanges $x_1^2$ and $x_3^2$; $s_{24}k_2k_4$
interchanges $x_1^2$ and $x_4^2$ -- more simply, consider the
symmetries of a genuinely black forked graph), so we need only
consider $\xi=\id$. Then $\bar e$ is as in~\eqref{projerepeated} and
and $K(\xi)=\{x_3^2,x_4^2\}$, and $\bar{\bar
e}=x_2^2+(x_2-x_4)^2+x_4^2$: $m(\xi)=x_4^2$ so $\gamma(\xi)=\delta(\xi)=2$.

Disconnected type: $\Span{x_1^2,x_2^2, (x_3-x_4)^2}$. In this case
$x_1^2$ and $x_2^2$ are interchanged by $k_1$ and $x_1^2$ and
$(x_3-x_4)^2$ are interchanged by $w'$, both of which
preserve~$\sigma$, so it is enough to consider $\xi=\id$. Then $\bar e$
is as in~\eqref{projerepeated} and $K(\xi)=\{x_2^2, (x_3-x_4)^2\}$,
and $\bar{\bar e}=2x_3^2+2x_4^2$: $m(\xi)=(x_3-x_4)^2$ so
$\gamma(\xi)=\delta(\xi)=4$.
\end{Proof}

\begin{cor}\label{mu3general}
If $C$ is a depth~$3$ curve with $\mu=3$ contained in the closure of
the image of the orbit of $\sigma$, and $\sigma$ is not of type
disconnected, suppose $H=aL-b D_4(n)-c E(n)$ is a divisor
on $\AVOR4(n)$ with $a-12 b/n\ge 0$, $b\ge 2c\ge 0$. Then $H.C\ge 0$.
\end{cor}
\begin{Proof}
Note first that $\sigma$ is not of type RT${}^*$, since the linear
forms whose squares span~$\sigma$ then span a space of dimension~$2$,
not~$3$. Since $\mu=3$ we have $\#K(\xi)=2$ always, by
Lemma~\ref{sizeofKxi}. According to Corollary~\ref{gammadeltacor} we
need to show that if $b\ge 2c$ then
$$
\sum_{\xi\in\gothS_3}\left(-\sfrac{b}{12}\delta(\xi)
-\sfrac{b}{4}\delta(\xi)
+\sfrac{b}{4}\gamma(\xi)\right)
+2b-c\ge 0.
$$
If $\sigma$ is of string type then
$\sum_{\xi\in\gothS_3}\delta(\xi)=12$ and
$\sum_{\xi\in\gothS_3}\gamma(\xi)=16$, so we get $2b-c\ge 0$; if
$\sigma$ is of type BF${}^*$ then
$\sum_{\xi\in\gothS_3}\delta(\xi)=\sum_{\xi\in\gothS_3}\gamma(\xi)=12$,
so we get $b-c\ge 0$.
\end{Proof}

If $\sigma$ is of disconnected type then
$\sum_{\xi\in\gothS_3}\delta(\xi)=\sum_{\xi\in\gothS_3}\gamma(\xi)=24$
and we get $-c$ which will not be positive. We need to deal with this
case separately. We do this by examining the contribution from
$D_{2,m(k)}(n)$. If $n\ge 3$ then the boundary component
$D_{2,m(k)}(n)$ of $\AVOR2(n)$ is isomorphic to the Shioda modular
surface (or universal elliptic curve) of level~$n$, which we call
$S(n)$: the projection to the modular curve is $s_{m(k)}\colon
D_{2,m(k)}(n) \to \AVOR1(n)=X(n)$.

\begin{lemma}\label{nonsection}
Let $N_2=-D_{2,m(k)}(n)|_{D_{2,m(k)}}(n)$ be the normal bundle of
$D_{2,m(k)}(n)$ in $\AVOR2(n)$. Then
$$
N_2=\sfrac{2}{n}s_{m(k)}^* L_{X(n)}+\sfrac{2}{n}\sum L_{ij}
$$
where the $L_{ij}$ are the sections of~$S(n)$.
\end{lemma}
\begin{Proof}
This was proved in~\cite[p. 271]{Hu}.
\end{Proof}
\begin{proposition}\label{mu3special}
Assume $C$ is a depth~$3$ curve with $\mu=3$ contained in the closure of
the image of the orbit of $\sigma$, and $\sigma$ is of disconnected
type. Suppose $H=aL-b D_4(n)-c E(n)$ is a divisor
on $\AVOR4(n)$ with $a-12 b/n\ge 0$, $b\ge 2c\ge 0$. Then $H.C\ge 0$.
\end{proposition}
\begin{Proof}
We may take $\sigma=\Span{x_1^2,x_2^2,(x_3-x_4)^2}$ and as in the
proof of Proposition~\ref{gammadeltavalues} we may assume
that~$\xi=\id$. Then the exceptional divisor is mapped by $r_\xi$ to
the line $\bar{\bar{E}}$ given by $\Span{x_3^2,x_4^2}$ and
$D_{2,m(k)}(n)$ intersects this line transversally. If
$r_\xi(C)\subset D_{2,m(k)}(n)$ is a curve that intersects~$\bar{\bar
E}$ then it does so at a singular point of a fibre of $s_{m(k)}$,
since the two components containing the image of~$E$ also meet
there. In particular, $r_\xi(C)$ cannot be a section of~$s_{m(k)}$.

We now have to remember that the term $A_2$ in
Proposition~\ref{HwithoutSxi} contains copies of the nef line bundle
$-D_{2,m(k)}(n)|_{D_{2,m(k)}}(n)$. In fact in view of the explicit description
of this line bundle given in Lemma~\ref{nonsection} we know that it is
represented by an effective $\QQ\,$-divisor, namely
$\frac{2}{n}(\sum L_{ij})+\frac{2}{n}s_{m(k)}^* L_{X(n)}=
\frac{2}{n}(\sum L_{ij})+s_{m(k)}^*\big(\frac{1}{6}X_\infty(n)\big)$, where
$X_\infty(n)$ is the set of cusps on the modular curve $X(n)=\AVOR1(n)$.
In particular, it follows that this line bundle has positive degree
on all curves that are not sections. In view of equation~\eqref{K2} and the
subsequent reasoning, it will be enough to prove that
$$
\bigg(-\sfrac{b}{12}\sum_{\xi\in\gothS_\mu}r_\xi^*D_{2,m(\xi)}(n) 
-\sfrac{b}{4}\sum_{\xi\in\gothS_\mu}r_\xi^*D_{2,m(\xi)}(n) 
-cE(n)_{\xi(1)\xi(2)}\bigg).C \geq 0,
$$
where $E(n)_{\xi(1)\xi(2)}$ denotes the union of those components of
$E(n)$ that map to $\bar{\bar E}$. This simplifies to
$$
\bigg(-\sfrac{b}{3}\sum_{\xi\in\gothS_\mu}r_\xi^*D_{2,m(\xi)}(n) 
-cE(n)_{\xi(1)\xi(2)}\bigg).C \geq 0.
$$
A priori the curve $C$ can meet several components of $E(n)_{\xi(1)\xi(2)}$
and each of these in several points. Let $P$ be some such point of
intersection and let $C_P$ be the curve germ defined by $C$ at~$P$.

Since $r_\xi(C)$ is not a section it follows that
$\frac{2}{n}(\sum L_{ij}).r_\xi(C)\geq 0$.
Now as we have said before the 
curve $r_\xi(C)$ must intersect a singular fibre of 
$D_{2,m(\xi)}(n)\cong S(n)$ in a singular point over some cusp. Denote
the fibre of $D_{2,m(\xi)}(n)$ over this cusp by $F_0$. Note that 
$D_{2,\hat{m}(\xi)}(n)$ intersects $D_{2,m(\xi)}(n)$ in two lines which
are contained in $F_0$ and which meet $\bar{\bar{E}}$.
Pulling back $-D_{2,\hat{m}(\xi)}(n)$ and using $\delta=4$, we obtain
$$
r_\xi^*(\sfrac{1}{6}F_0).C_P=\sfrac{1}{6}r_\xi^*(D_{2,\hat m(\xi)}).C_P 
\geq\sfrac{2}{3}E.C_P.
$$
But now the claim follows since
$$
-\sfrac{b}{3}\sum_{\xi\in\gothS_\mu}r_\xi^*D_{2,m(\xi)}(n).C_P\ge
\sum_{\xi\in\gothS_\mu}\sfrac{b}{18}r_\xi^*
F_0.C_P\ge \sfrac{4b}{3}E.C_P.
$$
So here $b \geq \sfrac{3}{4}c$ is enough.
\end{Proof}

\begin{lemma}\label{gammadeltaineq}
If $3 \le \mu \le 6$ and $\sigma$ is a cone such that the image of the
closure of its orbit contains a depth~$3$ curve, then $\gamma(\xi)\ge
\delta(\xi)= 2$ for all $\xi\in\gothS_\mu$.
\end{lemma}
\begin{Proof}
$\Sym_2\pr_{\xi(2)}$ maps $K(\xi)$ to $m(\xi)$. If $\bar e= \Sym\pr_{\xi(1)} e
={\bar l}_1^{\,2}+{\bar l}_2^{\,2}+ \bar l_3^{\,2} +\bar l_4^{\,2}$ then
$\gamma(\xi)=\psi_\gamma(\bar e)$ is equal to the number of the 
$\bar l_i^{\,2}$
that are not in~$K(\xi)$. But if $\bar l_i^{\,2}\in K(\xi)$, so
$\psi_\gamma(\bar l_i^{\,2})=0$, then
$\bar{\bar l}_i^{\,2}=(\pr_{\xi(2)}\bar l_i)^2=m(\xi)$ so
$\psi_\delta(\bar{\bar l}_i^{\,2})=0$, so $\psi_\delta(\bar{\bar e})\le
\psi_\gamma(\bar e)$. Hence $\gamma(\xi)\ge \delta(\xi)$.

But $\delta(\xi)$ is determined by $\xi(1)$ and $\xi(2)$ only, since
these determine $m(\xi)$ and hence $\psi_\delta$ as well as $\bar{\bar
e}$. Choose a non-RT${}^*$ face $\tau$ of $\sigma$ spanned by
$\xi(1)=l_1^2$, $\xi(2)=l_2^2$ and some other generator of
$\sigma$. These exist, since the third generator of an RT${}^*$ face
is in the linear span of $l_1$ and $l_2$, but the generators of
$\sigma$ span a rank~$3$ space. Moreover, such a non-RT${}^*$ face
cannot be of disconnected type, because any cone having a proper face
of disconnected type contains forms of rank~$4$. The choice of
$\xi(1)$, $\xi(2)$ and $\tau$ determines an element
$m(\xi_\tau)\in\MM_2$, namely
$m(\xi_\tau)=\Sym_2\pr_{\xi(1)\xi(2)}(\tau)$; and since $\tau$ is not
RT${}^*$ it is non-zero. Therefore $m(\xi_\tau)=m(\xi)$, since both of
them are the square of a generator of the same $1$-dimensional space
(namely, the projection of the linear space spanned by the generators
of~$\sigma$).

Now $\delta(\xi)$ is calculated exactly as one calculates $\delta$ for
the cone $\Span{\xi(1),\xi(2),\tau}$ with the generators in that
order. This is a non-RT${}^*$ cone with $\mu=3$ so by
Lemma~\ref{gammadeltavalues} that value of $\delta$ is equal to~$2$.
\end{Proof}

\begin{proposition}\label{mu45}
Suppose $C$ is a depth~$3$ curve with $\mu=4$ or $\mu=5$ and $H=aL-b
D_4(n)-c E(n)$ is a divisor on $\AVOR4(n)$ with $a-12 b/n\ge 0$, $b\ge
2c\ge 0$. Then $H.C\ge 0$.
\end{proposition}
\begin{Proof}
Since $\gamma(\xi)\ge\delta(\xi)$ we need only check that 
$$
\sum_{\xi\in\gothS_\mu}-\sfrac{b}{(\mu-1)\mu!}\delta_s(\xi)
+\sfrac{4b}{\mu-1}-c\ge 0.
$$
But $\delta(\xi)=2$ so 
$$
\sum_{\xi\in\gothS_\mu}-\sfrac{b}{(\mu-1)\mu!}\delta_s(\xi)
+\sfrac{4b}{\mu-1}-c= \sfrac{2b}{\mu-1}-c\ge 0.
$$
This is always fulfilled
for $\mu=4$ or $\mu=5$ and $b \geq 2c$.
\end{Proof}

\begin{proposition}\label{mu6}
Let $C$ be a depth~$3$ curve with $\mu=6$ and let $H=aL-b D_4-c E$ be a divisor
on $\AVOR4$ with $a\ge 0$, $a-12b\ge 0$ and $b\ge 2c\ge 0$. Then $H.C\ge 0$.
\end{proposition}
\begin{Proof}
In this case, by Lemma~\ref{mu6face}, we can always work with the cone 
$\sigma = \Span{x_1^2,x_3^2,x_4^2,(x_1-x_3)^2,(x_1-x_4)^2,(x_3-x_4)^2}$.
However, we cannot restrict ourselves to $\xi=\id$.

There are two essentially different cases, depending on whether the
edges $\xi(1)$ and $\xi(2)$ in the bicoloured
graphs given in the proof of Lemma~\ref{mu6face} (Figure~4) are
opposites ($3$ cases) or not ($12$ cases): see
Remark~\ref{opposites}. 
A representative of the 
non-opposite case is $(\xi(1),\xi(2))=(x_1^2,x_3^2)$, of the other 
$(\xi(1),\xi(2))=(x_1^2,(x_3-x_4)^2)$. The geometry of these cases is as
follows. In the non-opposite case we obtain
$$
\bar{e}=(x_2-x_3)^2+(x_2-x_4)^2+x_3^2+x_4^2,
$$
$$
\bar{\bar{e}}=x_2^2+(x_2-x_4)^2+x_4^2,
$$
$$
m(\xi)=x_4^2.
$$
In this case the exceptional divisor is mapped to a point $\bar{\bar{E}}$
which is a singular point of a singular fibre of $D_{2,m(\xi)}(n)$
and if $C$ meets this divisor then $r_\xi(C)$ cannot be a section.
In the opposite case we have
$$
\bar{e}=(x_2-x_3)^2+(x_2-x_4)^2+x_3^2+x_4^2,
$$
$$
\bar{\bar{e}}=2(x_2-x_4)^2+2x_4^2,
$$
$$
m(\xi)=x_4^2.
$$
In this case $\bar{\bar{E}}$ is a line which is contained in
$D_{2,m(\xi)}(n)$ and we cannot a priori exclude that $r_\xi(C)$ is a section.

It is straightforward to check by hand that this distinction
coincides with the distinction in terms of opposite and non-opposite
edges above, and to verify directly that there are $12$~non-opposite
cases.

We shall now argue as in the case $\mu=3$ of disconnected type
(Proposition~\ref{mu3special}), but taking into account only the
contribution from the non-opposite cases. The other contributions are
non-negative. Using the same notation as in the proof of
Proposition~\ref{mu3special} and using that $\delta=2$ gives us
$-r_\xi^*D_{2,m(\xi)}(n).C_P \geq
r_\xi^*(\sfrac{1}{6}F_0).C_P\geq\sfrac{1}{3}E.C_P$.
In view of formula~\eqref{K3} it will be enough to prove that
$$
\bigg(-\sfrac{b}{5\cdot 6!}\sum_{\xi\in\gothS'_\mu}r_\xi^*D_{2,m(\xi)}(n)
-\sfrac{b}{10\cdot 5!}\sum_{\xi\in\gothS'_\mu}r_\xi^*D_{2,m(\xi)}(n)
+(\sfrac{2b}{5}-c)E(n)\bigg).C_P \geq 0,
$$
where $\gothS'_{\mu}$ is the set of $\xi$ giving rise to the
non-opposite case, so $\#\gothS'_\mu=\frac{4}{5}6!$. This simplifies to
$$
-\sfrac{4b}{5\cdot 6!}\sum_{\xi\in\gothS'_\mu}r_\xi^*D_{2,m(\xi)}(n).C_P
+(\sfrac{2}{5}b-c)E.C_P \geq 0.
$$
This leads to $b \geq \frac{75}{46}c$, and $\frac{75}{46}<2$ so we are done.
\end{Proof}

\begin{proposition}\label{nefondepth3}
Let $C$ be a depth~$3$ curve and let $H=aL-b D_4-c E$ be a divisor
on $\AVOR4$ with $a\ge 0$, $a-12b\ge 0$ and $b\ge 2c\ge 0$. Then $H.C\ge 0$.
\end{proposition}
\begin{Proof}
This follows from the cases dealt with above, in
Propositions~\ref{mu3general}, \ref{mu3special}, \ref{mu45} and \ref{mu6}.
\end{Proof}

Theorem~\ref{mainthm} now follows from Propositions~\ref{nefondepth4},
\ref{nefondepth0}, \ref{nefondepth1}, \ref{nefondepth22},
\ref{nefondepth23} and \ref{nefondepth3}.


\bibliographystyle{alpha}

\bigskip

\begin{tabular}{ll}
K. Hulek & G.K. Sankaran\\
Institut f\"ur Mathematik & Department of\\
Universit\"at Hannover & \ \ Mathematical Sciences\\
Postfach 6009 & University of Bath \\
D 30060 Hannover & Bath BA2 7AY\\
Germany & England\\
\\
{\tt hulek@math.uni-hannover.de} & {\tt gks@maths.bath.ac.uk}
\end{tabular}

\end{document}